\RequirePackage{snapshot}
\documentclass{scrartcl}

\usepackage[utf8]{inputenc}
\usepackage[english]{babel}
\usepackage[margin=2cm]{geometry}

\usepackage[backref,colorlinks=true, urlcolor=black, linkcolor=black, citecolor=black]{hyperref}

\usepackage{iftex}

\ifluatex
  \usepackage{mathtools}
  \usepackage{unicode-math}
\else
  \usepackage{mathtools}
  \usepackage{amsmath, amssymb}
  \usepackage{stix2}
  \usepackage[T1]{fontenc}
\fi

\usepackage{amsthm}
\usepackage{commath}
\usepackage[numbers]{natbib}
\usepackage{graphicx}
\usepackage[dvipsnames]{xcolor}
\usepackage{booktabs}
\usepackage{placeins} 
\usepackage{subcaption}
\usepackage{algpseudocode}
\usepackage{algorithm}
\usepackage{framed}
\usepackage[normalem]{ulem}
\usepackage{cancel}
\usepackage{enumerate}
\usepackage{cleveref}
\usepackage{wrapfig}
\usepackage{tikz}\usetikzlibrary{positioning, shapes.geometric, shapes, arrows.meta, calc}
\usepackage{pgfplots}
\usepackage{ifthen}
\pgfplotsset{width=10cm,compat=1.9}
\usepackage{listofitems} 
\usepackage[outline]{contour} 
\contourlength{1.4pt}

\tikzset{>=latex} 

\usepackage{xcolor}
\colorlet{myred}{red!80!black}
\colorlet{myblue}{blue!80!black}
\colorlet{mygreen}{green!60!black}
\colorlet{myorange}{orange!70!red!60!black}
\colorlet{mydarkred}{red!30!black}
\colorlet{mydarkblue}{blue!40!black}
\colorlet{mydarkgreen}{green!30!black}
\tikzset{
  >=latex, 
  node/.style={thick,circle,draw=myblue,minimum size=22,inner sep=0.5,outer sep=0.6},
  node in/.style={node,green!20!black,draw=mygreen!30!black,fill=mygreen!25},
  node hidden/.style={node,blue!20!black,draw=myblue!30!black,fill=myblue!20},
  node convol/.style={node,orange!20!black,draw=myorange!30!black,fill=myorange!20},
  node out/.style={node,red!20!black,draw=myred!30!black,fill=myred!20},
  connect/.style={thick,mydarkblue}, 
  connect arrow/.style={-{Latex[length=4,width=3.5]},thick,mydarkblue,shorten <=0.5,shorten >=1},
  node 1/.style={node in}, 
  node 2/.style={node hidden},
  node 3/.style={node out}
}
\def\nstyle{int(\lay<\Nnodlen?min(2,\lay):3)} 
\def\n2style{int(\lay<\Nnodlen+4?min(2,\lay-4):3)} 

\newcommand{\R}{\mathbb{R}}
\newcommand{\N}{\mathbb{N}}

\newcommand{\udef}{\coloneqq}

\theoremstyle{definition}

\newcommand{\be}[1]{
\begin{equation}
\expandafter\label{eq:#1}
}
\newcommand{\ben}{
\begin{equation}
\expandafter
}
\newcommand{\ee}{\end{equation}}

\newcommand{\bg}[1]{
\begin{gather}\nonumber
\expandafter\label{eq:#1}
}
\newcommand{\eg}{\end{gather}}

\newcommand{\bfig}[1]{
\begin{figure}
\expandafter\label{fig:#1}
}
\newcommand{\efig}{\end{figure}}

\def\({ \left( }
\def\){ \right) }

\newcommand{\x}{x}
\newcommand{\y}{y}

\renewcommand{\c}{x}

\renewcommand{\c}{c}

\newcommand{\diverg}{\nabla\cdot}

\newcommand{\grad}{\nabla}

\newcount\colveccount
\newcommand*\colvec[1]{
        \global\colveccount#1
        \begin{pmatrix}
        \colvecnext
}
\def\colvecnext#1{
        #1
        \global\advance\colveccount-1
        \ifnum\colveccount>0
                \\
                \expandafter\colvecnext
        \else
                \end{pmatrix}
        \fi
}

\newcount\rowveccount
\newcommand*\rowvec[1]{
        \global\rowveccount#1
        \begin{pmatrix}
        \rowvecnext
}
\def\rowvecnext#1{
        #1
        \global\advance\rowveccount-1
        \ifnum\rowveccount>0
                &
                \expandafter\rowvecnext
        \else
                \end{pmatrix}
        \fi
}

\graphicspath{{./Figures/}}

\title{
Inverse Physics-Informed Neural Networks for transport models in porous materials
}
\author{Marco Berardi
\footnote{Consiglio Nazionale delle Ricerche, Istituto di Ricerca Sulle Acque,  Bari, Italy, \url{marco.berardi@cnr.it}}
,
Fabio V. Difonzo
\footnote{Department of Engineering,
LUM University Giuseppe Degennaro,
S.S. 100 km 18, 70010 Casamassima (BA), Italy, \url{difonzo@lum.it}}
\footnote{Consiglio Nazionale delle Ricerche,  Istituto per le Applicazioni del Calcolo "Mauro Picone", Bari, Italy, \url{fabiovito.difonzo@cnr.it}}
,
Matteo Icardi
\footnote{School of Mathematical Sciences, University of Nottingham, Nottingham, UK, \url{matteo.icardi@nottingham.ac.uk}
}
\footnote{Universit\'a degli Studi di Bari Aldo Moro, Dipartimento di Matematica, Bari, Italy}
}
\date{\today}
\begin{document}

\maketitle

\section*{Abstract}
Physics-Informed Neural Networks (PINN) are a machine learning tool that can be used to solve direct and inverse problems related to models described by Partial Differential Equations by including in the cost function to minimise during training the residual of the differential operator. This paper proposes an adaptive inverse PINN applied to different transport models, from diffusion to advection-diffusion-reaction, and mobile-immobile transport models for porous materials. Once a suitable PINN is established to solve the forward problem, the transport parameters are added as trainable parameters and the reference data is added to the cost function. We find that, for the inverse problem to converge to the correct solution, the different components of the loss function (data misfit, initial conditions, boundary conditions and residual of the transport equation) need to be weighted adaptively as a function of the training iteration (epoch). Similarly, gradients of trainable parameters are scaled at each epoch accordingly.  Several examples are presented for different test cases to support our PINN architecture and its scalability and robustness.



\section{Introduction}


In recent years, Physics-Informed Neural Networks (\textsc{PINNs}) (see \citep{karniadakis2021physics} for a recent review) have attracted significant attention in mathematical modelling due to their ability to address direct and inverse problems governed by differential equations. This tool elegantly integrates the principles of physics-based differential models with the adaptability of neural network architectures.

On the one hand, \textsc{PINNs} offer a powerful framework for solving direct problems: those concerned with computing the solution of complex partial differential equations (PDEs) with defined initial and boundary conditions. Classical numerical techniques, such as finite difference methods, finite elements, virtual element schemes, and spectral methods, are well established. However, these methods may encounter difficulties when addressing problems characterised by high non-linearities, high dimensionality, or uncertainties in parameters or boundary conditions. To overcome these limitations, data-driven approaches have been explored. For example, \citet{Zhou2023} explored using deep neural networks with specialised activation functions for solving high-dimensional nonlinear wave equations. \citet{SUKUMAR2022} proposed a geometry-aware \textsc{PINN} method specifically designed to enforce boundary conditions in complex domains. Extensions to integral equations also show promise, as demonstrated by \citet{Zunino}, where orthogonal decomposition is combined with neural networks.

On the other hand, \textsc{PINNs} can play a crucial role in solving inverse problems as surrogate models coupled with standard parameter estimation techniques (Bayesian or deterministic), or as stand-alone tools. In the first case, the independent variables and the physical parameters are treated as input to the neural network, which is trained for a wide range of physical parameters. In the second case, the physical parameters are treated as trainable neural network parameters. Although not directly appearing in the neural network, they are present in the loss function through the equation residual. In this work we will consider only this second approach which has been successfully applied to a wide range of problems, from groundwater flow and contaminant transport \citep{DLP2024,CUOMO2023106}, to heat transfer in porous media \citep{Yang_Karniadakis_JPC_2021}, to the identification of unknown PDE structures \citep{GAO2022}. We will refer to this approach as \emph{inverse PINN}.
 Inverse problems are inherently more complex than their direct counterparts due to their potential ill-posed nature, where multiple solutions might exist or none at all. Additional challenges arise from data-scarce regimes, irregular geometries (e.g., \citep{GAO2022}), missing data, or uncertainties inherent to the model. Advanced \textsc{PINN} techniques have been developed to address these issues. \citet{Yang_Karniadakis_JPC_2021} introduced a Bayesian \textsc{PINN} framework for both inverse and forward models. \citet{Gusmao_Medford_2024} framed \textsc{PINNs} in terms of maximum-likelihood estimators, enabling error propagation and removing a hyperparameter. Finally, \citet{DLP2024} employs a serialised \textsc{PINN} approach to determine kernel functions in peridynamic models, arising from non-local formulation of continuum mechanics (e.g. \cite{LPcheby,BDFP}).
 
Data-driven and machine-learning approaches have found promising applications particularly in material modelling \citep{jeong2024data,guo2022machine,upadhyay2024physics} and transport processes in porous media \citep{amini2022physics,marcato2023reconciling,santos2020poreflow}. We refer to \citep{d2022machine} and references therein for a recent review. More specifically, porous media are present in almost all aspects of engineering, manufacturing, and physical sciences. Porous media effective parameters (i.e. parameters controlling the emerging macroscopic dynamic of multi-scale materials) often have to be found by solving inverse problems. These include, for example, permeability, dispersivity and effective reactivity, which are crucial for agronomy, soil science and hydrological applications (e.g. \citep{BerardiDifonzoGuglielmi2023,Celia_et_al,KUMAR_POP_RADU,Marinoschi_book}). 
In fact, direct measures of these parameters, such as hydraulic conductivity and porosity, can be time-consuming, expensive and difficult to spatialize: for this reason, inferring these parameters from more easily measurable quantities, for instance representing the state variables in a process governed by differential systems, can be helpful to a more significant assessment (see for instance \cite{DiLena_Berardi_Masciale_Portoghese_HYP,DuChateau_SIAM_1997,Rossi_et_al_2015}). A proper estimation of parameters is mandatory for correctly forecasting the dynamics of significant processes in porous media; for instance, the saturated hydraulic conductivity plays a crucial role in the dynamics of soil moisture content, and this parameter can vary up to many orders of magnitude, as in \cite{Berardi_Difonzo_Vurro_Lopez_ADWR_2018}.
Similar problems arise also in biology and medicine (tissues,  bones, circulation network) and engineering (porous electrodes in batteries, concrete and building materials and materials design, e.g \citep{Wein_et_al_CNSNS_2019,Frittelli_Madzvamuse_Sgura_2023}). Due to their multi-scale structure and heterogeneity and the availability of sparse heterogeneous data, data-driven models like PINNs quickly gained popularity in these areas. Recent works also investigated the use of PINNs for the identification of the effective parameters \citep{he2020physics,wang2023physics,he2021physics,bandaiGhezzehei2022}.\\
{Different methods for handling inverse problems have been proposed and used in last decades. For instance, Kalman filters have been widely used for parameter estimation, as in \cite{Crestani_et_al_hess_2013,Medina_et_al_hess_2014}, due to their ease of implementation. Bayesian methods for estimating parameters have also been introduced in \cite{Kennedy_OHagan_2002}, significantly reviewed in \cite{Viana_Subramaniyan_2021}, and exploited in a porous media context, for instance, in \cite{Xu_et_al_WRR_2017}. 
Unlike Kalman filters and Bayesian methods, which require numerous forward model evaluations and a-priori assumptions about the distribution of the parameters
 (see, for a thorough review, \cite{Humpherys_et_al_SIAM_review_2012}), PINNs are based on the minimisation of a deterministic loss function, and are differentiable with respect to the parameters to be estimated. This makes them particularly suitable for inverse problems.\\
Nevertheless, comparing the efficiency of such different methods for parameters estimation is beyond the scopes of this paper, which aims at proposing a novel approach for inverse PINN.}

This work proposes an adaptive inverse \textsc{PINN} architecture for solving advection-diffusion-reaction models. We consider a number of transport models, from simple diffusion to more complex advection-diffusion-reaction equations with the key parameters being dispersivity, effective velocity (a measure of permeability), and reaction constants. The main novelty of the proposed approach is the adaptive scaling of the loss function components and gradients of the trainable parameters. This adaptive scaling is crucial for the convergence of the inverse problem, as it ensures that the different components of the loss function are balanced and that the gradients of the trainable parameters are scaled appropriately. We demonstrate the effectiveness of the proposed approach through a series of numerical experiments, showing that the adaptive inverse \textsc{PINN} architecture is scalable, robust, and efficient for solving a wide range of transport models.
 In \cref{sec:models}, we introduce the mathematical models we consider, and in \cref{sec:PINN}, we present the \textsc{PINN} architecture we use. In \cref{sec:results}, we present the results of our numerical experiments, and in \cref{sec:conclusions}, we draw our conclusions. All codes and data used in this work are published and freely available \citep{icardi_2024_12743929}.

\section{Mathematical models}
\label{sec:models}

\subsection*{Heat equation}
\label{sec:heat}
The simplest model we consider is a pure diffusion (heat) equation:
\begin{align}
        \label{eq:heat}
        u_t - D\Delta u &=  0\,,\qquad x\in\Omega\,,\; t\in[0,T]\,,
\end{align}
and its non-linear extension:
\begin{align}
        \label{eq:heatnl}
        u_t - \diverg(\mathcal{D}(u)\grad u) &=  0\,,\qquad x\in\Omega\,,\; t\in[0,T]\,,  
\end{align}
where $D$ is a constant diffusion coefficient, and $\mathcal{D}(u)$ is a non-linear diffusion coefficient.

\subsection*{Advection-Diffusion-Reaction equation}
\label{sec:adr}
We consider the following advection-diffusion-reaction equation with non-linear reaction term:
\begin{align}
        \label{eq:adr}
        \beta u_t + \diverg\left(V u - D\grad u\right) &=  \sigma(u)\,,\qquad x\in\Omega\,,\; t\in[0,T]\,,
\end{align}
with appropriate initial and boundary conditions. The spatial operator is a linear advection diffusion operator, with velocity $V$ and dispersion/diffusion coefficient $D$. The time-derivative involves a porosity term $\beta$.
The reaction term $\sigma(u)$ is a non-linear function of the concentration $u$, which can be used to model a wide range of physical and chemical processes.

\subsection*{Mobile-Immobile model}
\label{sec:mi}

\begin{figure}[htbp]
        \centering
        \includegraphics[width=0.6\textwidth]{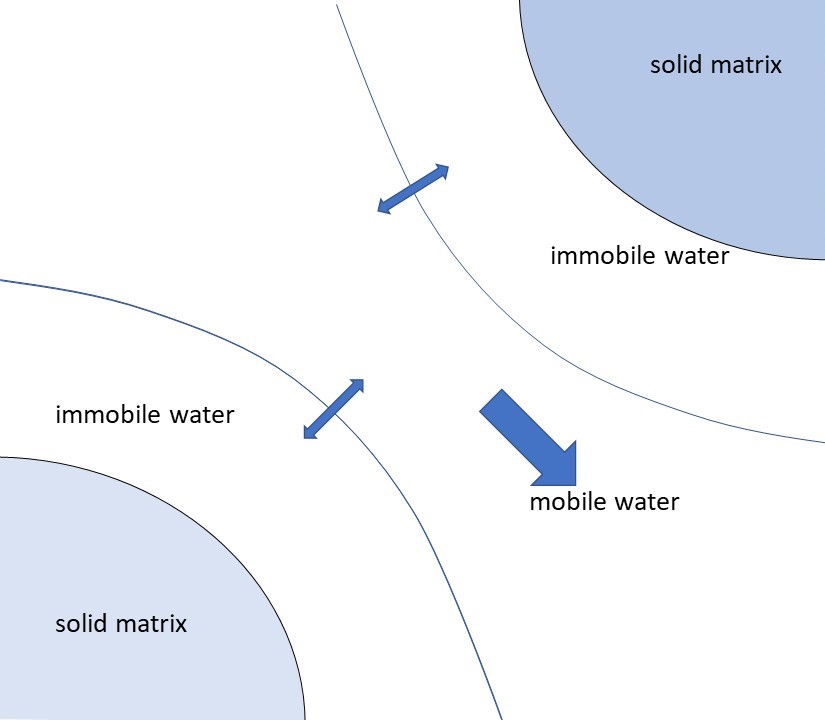}
        \caption{A graphical representation of the mobile-immobile model for the transport of solutes in porous media: as in \cite{DeSmedt_Wierenga_1979}, the mobile region is the primary zone of water and solute transport; the other region is termed   the immobile zone, because the soil water in this zone is stagnant relative  to the water in the mobile zone.}
        \label{fig:mobile_immobile}
\end{figure}

The Mobile-Immobile model \citep{lapidus1952mathematics} describes the transport of solutes in porous media and is based on the assumption that the solute is partitioned between a mobile and an immobile phase (see  \cref{fig:mobile_immobile}). The mobile phase is assumed to be in equilibrium with the immobile phase, and the transfer of solute between the two phases is described by a first-order rate equation. Historically, this model has been developed after observing an anomalous behaviour of solute breakthrough curves measured at field and laboratory scale (see, for instance \cite{Masciopinto_Passarella_2018}). The model has been applied to a wide range of problems, including the transport of contaminants in groundwater, the transport of nutrients in soil, and the transport of solutes in fractured rock. For a recent review, derivation and extensions of the model we refer to \citep{municchi2020generalized,dentz2018mechanisms}.
The model can be written as follows:
\begin{align}
        \label{eq:mi}
        \beta_0 u_t + \diverg\left(V u - D\grad u\right) &= \lambda(v-u)\,,\qquad x\in\Omega\,,\; t\in[0,T]\,,
        \\ \nonumber
        \beta_1 v_t &= -\lambda(v-u)\,,\qquad x\in\Omega\,,\; t\in[0,T]\,,
\end{align}
where $u$ and $v$ are the concentrations of the solute in the mobile and immobile water phases, respectively,
$\beta_0$ and $\beta_1$ are the mobile and immobile porosities (volume fractions), and $\lambda$ is the transfer coefficient describing the rate of transfer of solute between the mobile and immobile phases. The transfer coefficient is assumed to be constant in time and space.

\section{Physics-Informed Neural Network}
\label{sec:PINN}

\begin{figure}
        \centering
        \begin{tikzpicture}[x=1.7cm,y=1.4cm]
          \message{^^JNeural network, shifted}
          \readlist\Nnod{2,4,4,4,1} 
          \def\yshift{0.5} 
          
          \message{^^J  Layer}
          \foreachitem \N \in \Nnod{ 
            \def\lay{\Ncnt} 
            \pgfmathsetmacro\prev{int(\Ncnt-1)} 
            \message{\lay,}
            \foreach \i [evaluate={\c=int(\i==\N); \y=\N/2-\i-\c*\yshift;
                         \x=\lay; \n=\nstyle;}] in {1,...,\N}{ 
              \node[node \n] (N\lay-\i) at (\x,\y) {}
              ;
              \ifnum\lay>1 
              \ifnum\lay<4
                \foreach \j in {1,...,\Nnod[\prev]}{ 
                  \draw[connect arrow] (N\prev-\j) -- (N\lay-\i);
                }
              \fi
              \fi 
              \ifnum\lay>4
                \foreach \j in {1,...,\Nnod[\prev]}{ 
                  \draw[connect arrow] (N\prev-\j) -- (N\lay-\i);
                }
              \fi
              \ifnum\lay=4
                \foreach \j in {1,...,\Nnod[\prev]}{ 
                  \draw[dotted] (N\prev-\j) -- (N\lay-\i);
                }
              \fi
              
            }
            \ifnum\lay>1 \ifnum\lay<5
            \path (N\lay-\N) --++ (0,1+\yshift) node[midway,scale=1.5] {$\vdots$};
            \fi\fi
          }
          
          \node[left=0.6,align=center,mygreen!60!black] at (N1-1.0) {space\\[-0.2em]input\\[-0.2em]$x$};
          \node[left=0.6,align=center,mygreen!60!black] at (N1-2.0) {time\\[-0.2em]input\\[-0.2em]$t$};
          \node[above=0.6,align=center,myblue!60!black] at (N3-1.-90) {\texttt{tanh}\\[-0.2em] hidden layers\\[-0.2em]$\theta=\{W_1,b_1,\dots,W_l,b_l,\dots,W_L,b_L\}$
          };
          \node[right=0.2,align=center,myred!60!black] (outputlayer) at (N\Nnodlen-1.0) {output\\[-0.2em]layer\\[-0.2em]$u_{NN}(x,t;\theta_0)$};

          \node[draw, rectangle, rounded corners=5pt, align=center, below=1cm] (eqbox) at (N3-4.-90) {$\displaystyle\arg\min_{u,\theta_0}\mathcal{L}(u;\theta_0)$};
          \draw[->, thick, bend left] (outputlayer.south) to (eqbox.east);

          \node[below left=0.5cm of eqbox, align=center, myred!60!black] (direct) {$u$\\Direct problem};
          \node[below right=0.5cm of eqbox, align=center, mydarkred!60!black] (inverse) {$\theta_0$\\Inverse problem};

          \draw[->,thick,myblue] (eqbox.north) -- (N3-4.south);

          \draw[->, thick] (eqbox.south) -- (direct.north);
          \draw[->, thick] (eqbox.south) -- (inverse.north);
          
        \end{tikzpicture}
        \caption{PINN structure used in this work, with $L$ layers, $m_l$ neurons per layer, and hyperbolic tangent activation function in the hidden layers. The set $\theta_0\subseteq\{V,D,\lambda\}$ contains the trainable parameters relative to the considered physical law.}

        \label{fig:NNstructure}
        \end{figure}
      
In this paper, we will consider a Feed-Forward fully connected Neural Network (FF-DNN), also called Multi-Layer Perceptron (MLP) (see \citep{bengio2003,Cybenko1989} and references therein). \\
In a PINN, the solution space is approximated through a combination of activation functions, acting on all the hidden layers, with the independent variable used as the network input. Letting $(x,t)\in\R^{d+1}$ be the input of the NN, in a Feed-Forward network, each layer feeds the next one through a nested transformation so that it can be expressed, letting $L$ be the number of layers, as 
\begin{equation}\label{eq:NN}
\begin{aligned}
z_0 &= (x,t), \\
z_l &= \rho_l\left(\Lambda_l(z_{l-1})\right)\,,\qquad l=1,\ldots,L, \\
\Lambda_l(z_{l-1})&\udef W_lz_{l-1}+b_l 
\end{aligned}
\end{equation}
where, for each layer $l=1,\ldots,L$, $\rho_l:\R^{m_{l}}\to\R^{m_l}$ is the activation function, which operates componentwise, $W_l\in\R^{m_l\times m_{l-1}}$ is the weight matrix and $b_l\in{\R^{m_{l}}}$ is the bias vector. Thus, the output $z_L\in\R^m$ of a FF-NN can be expressed as a single function of the input vector $x$, defined as the composition of all the layers above in the following way: 
\[
u_{NN}(x,t;\theta)=z_L\udef(\rho_L\circ\Lambda_L\circ\ldots\circ\rho_1\circ\Lambda_1)(x,t).
\]
We denote the training parameters set as $\theta= \{W_l,b_l\}_{l=1}^L$. \\
In \cref{fig:NNstructure} we show a schematic representation of the structure of a PINN, where the input layer is composed of two neurons, one for the spatial variable and one for the time variable, and the output layer is composed of a single neuron, resulting in the approximation $u_{NN}(x,t)$. The hidden layers are composed of the same number of neurons $m_l$, $l=1,\ldots,L-1$, and the activation function used in the hidden layers is the hyperbolic tangent function, while the output layer uses the identity function. We have, therefore, $m_0=d+1$, $m_L=1$, and $m_l=m$ for $l=1,\ldots,L-1$. The activation function $\rho_l$ is the hyperbolic tangent function for each hidden layer $l=1,\ldots,L-1$, and the identity function for $l=L$. We will limit to one-dimensional spatial domain, i.e., $d=1$. The specific values for $m$ and $L$ will differ for each case in Section \ref{sec:results} and will be therefore given thereafter.

The aim of a PINN is to minimise  a suitable objective function called \emph{loss function} that includes not only the the data but also the physics of the problem. The minimisation is performed with respect to all the trainable parameters $\theta$, through a Stochastic Gradient Descent method.
Given a general spatio-temporal differential operator $\mathcal{D}(u;\theta_0)=0$, where $\mathcal{D}$ represents the differential operator acting on the unknown function $u\in V(\R^d)$, with physical parameters $\theta_0\in\R^s$, the loss function used by a PINN is given by
\begin{equation}\label{eq:lossPINN}
\mathcal{L}(u;\theta_0)\udef\sum_{i=1}^{M}\left(\|u(x^*_i,t^*_i)-u^*_i\|^2+\|\mathcal{D}(u(x^*_i,t^*_i);\theta_0)\|^2\right),
\end{equation}
where $u^*_i$ is the unknown function measured at point $(x^*_i,t^*_i)$ inside the domain or on the boundary. \\
The set $X^*=\{(x^*_i,t^*_i)\}_{i=1}^M$ is the set of training points, and $M$ is the number of training points. We highlight here that, since we are going to compare the PINN solution to synthetic data, we select collocation points coincident with training points. We recall that training points are used to teach the network to fit the known solution in the data-driven regions of the problem space, whereas collocation points are used to ensure that the solution, provided by the neural network, respects the physical law modelled by the differential equation considered (for further details we refer to \cite{lau2024pinnacle}). Therefore, in general training points and collocation points could be different collocation points are a subset of training points but, here, they will be given by the same set of points. \\
The chosen norm $\|\cdot\|$ (it may be different for each term in the loss function) depends on the functional space and the specific problem. Selecting a correct norm (to avoid overfitting) for the loss function evaluation is an important problem in PINN, and recently, in \citep{TAYLOR2023}, the authors have proposed spectral techniques based on Fourier residual method to overcome computational and accuracy issues. The first term in the right-hand side of \cref{eq:lossPINN} is referred to as data fitting loss and could possibly handle both initial and boundary conditions, while the second term is referred to as residual loss, which is responsible for making the NN informed by the physics of the problem.
The derivatives inside $\mathcal{D}$ in space, time and in the parameter space are usually performed using \texttt{autodiff} (Automatic Differentiation algorithm, see \citep{BaydinEtAl2017,BoltePauwels2020}).
Using the NN to approximate $u$ in the loss function \cref{eq:lossPINN} allows us to solve the PDE by minimising the loss function with respect to the parameters $\theta$ of the NN. If $u\approx u_{NN}(\theta)$, then the minimisation problem can be written as
\begin{equation}\label{eq:inverseProblem}
\theta^\dagger=\arg\min_{\theta}\mathcal{L}(u_{NN}(\theta);\theta_0).
\end{equation}
For a more detailed discussion on the PINN structure and the loss function, we refer to \citep{RAISSI2019}, \citep{Yang_Karniadakis_JPC_2021} and to the review in \citep{Cuomo2022}.

\subsection{Inverse PINN}
\label{sec:inversePINN}
Inverse PINNs are a type of Neural Network specifically designed to determine constitutive parameters or problem-related functions that appear in the PDE one must solve. However, due to the limited amount of data relative to exact solutions, or of available measurements of the physical problem described by the PDE underlying the PINN, the inverse problem \eqref{eq:inverseProblem} could likely be ill-posed, and thus particular care has to be put into the training strategy during the optimization process (see, e.g., \cite{Arridge_Maass_Öktem_Schönlieb_2019}). In particular, different contributions in the loss functions \cref{eq:lossPINN} could conflict with each other, providing an unbalanced gradient back-propagation during the training, which would result in a troublesome convergence process \citep{XU2023}. Thus, several strategies have been recently developed to cope with these issues: among the others, one could resort to \emph{GradNorm} \citep{HAGHIGHAT2021} to dynamically tune gradient magnitudes to balance learning tasks; to \emph{PCGrad} \citep{yu2020gradient} to project each gradient on the tangent plane to all the other conflicting gradients to mitigate such destructive interference; to \emph{Multi-Objective Optimization} \citep{sener2018multi}; to \emph{Self-Adaptive PINNs} \citep{mcclenny2022selfadaptive}, where each training point is weighed individually, so to penalize more points in difficult regions of the domain. 

Using the notation introduced in the previous section, the inverse PINN minimisation now takes into account also the physical parameters and can be written as
\begin{equation}
(\theta^\dagger,\theta_0^\dagger)= \left[\arg\min_{\theta,\theta_0}\mathcal{L}(u_{NN}(\theta);\theta_0)+\iota\|\theta_0-\theta_0^*\|^2\right],
\end{equation}
where $\iota$ is a regularisation parameter. The second term in the right-hand side of the equation is the regularisation term, which is used to prevent overfitting and to ensure that the physical parameters $\theta_0$ are close to the some reference parameters $\theta_0^*$.  If otherwise stated, in the following we will consider $\iota=0$.

With reference to the mathematical models introduced in \cref{sec:models}, we will consider the following physical parameters to be estimated: the diffusion coefficient $D$ in the heat equation \eqref{eq:heat}, the velocity $V$ and the dispersion coefficient $D$ in the advection-diffusion-reaction equation \eqref{eq:adr}, and the transfer coefficient $\lambda$ in the mobile-immobile model \eqref{eq:mi}. The physical parameters will be considered as trainable parameters in the NN, and the reference data will be added to the loss function \cref{eq:lossPINN}.

\subsection{Adaptive inverse PINN}
\label{sec:adaptivePINN}

To ensure the convergence of the inverse PINN, we redefine a weighted loss function as
\begin{equation}
\mathcal{L}(u_{NN}(\theta);\theta_0)=
\sum_{i=1}^{M}
\left(
        \omega_i^k\|u(x^*_i,t^*_i)-u^*_i\|^2
        +
        \omega_{\mathcal{D}}^k \|\mathcal{D}(u(x^*_i,t^*_i);\theta_0)\|^2
\right),
\end{equation}
where $\omega_i^k$ are weight factors that depends on the training iteration $i$. The weights are updated at each iteration to ensure that the different components of the loss function are balanced. The weights are updated using the following formula:
\begin{align}\label{eq:omega}
        \omega_i^{k} &= \frac{\hat{\omega}_i^k}{\sum_{j=1}^{M}\hat{\omega}_j^k + \hat{\omega}_{\mathcal{D}}^k}\,,\qquad i=1,\ldots,M\,,
        \\
        \omega_{\mathcal{D}}^{k} &= \frac{\hat{\omega}_{\mathcal{D}}^k}{\sum_{j=1}^{M}\hat{\omega}_j^k + \hat{\omega}_{\mathcal{D}}^k}\,,
        \\
        \hat{\omega}_i^{k} &=
        \begin{cases}
                \eta_{BC} & \text{if } x^*_i \in \partial\Omega \\
                \eta_{IC} & \text{if } t^*_i = 0 \\
                \nu(k) \eta_{u} & \text{if } (x^*_i,t^*_i) \text{ is a collocation point} \\
                0 & \text{otherwise}
        \end{cases}
        \,,\qquad i=1,\ldots,M\,,
        \\
        \hat{\omega}_{\mathcal{D}}^{k} &= 1,
\end{align}
where $\eta_{BC}$, $\eta_{IC}$, and $\eta_{u}$ are the weights for the boundary conditions, initial conditions, and collocation points, respectively. The function
$\nu_k$ is an increasing function of the epoch $k$, such that $\nu_k=0$ and $\nu_k\to 1$ as $k\to\infty$. This allows the PINN to be trained initially solely by the PDE residual. 
In the following, we will consider
\begin{equation}\label{eq:nu_k}
\nu(k) = \frac{{\tanh\left(10\left(\frac{{k-K/2-K_0}}{K}\right)\right)+1}}{2}
\,,\qquad k=1,\ldots,K\,,
\end{equation}
where $K$ is the total number of epochs, and $K_0$ is a threshold epoch before the weights are updated more significantly; see \cref{fig:nu_k} for a typical graph of a function of this kind. 

\begin{figure}
\centering
\scalebox{0.5}{
\begin{tikzpicture}[
  declare function={
    func(\x)= (tanh((\x-3500)/500)+1)/2;
  }
]
\begin{axis}[
    width = \linewidth,
    axis lines = left,
    xlabel = $k$,
    ylabel = {$\nu(k)$},
    ymin = 0,
]

\addplot[
    domain=1:5000, 
    samples=200, 
    color=black,
] {func(x)};
\end{axis}
\end{tikzpicture}
}
\caption{Qualitative behaviour of $\nu_k$ in \eqref{eq:nu_k} for $K=5000$ and $K_0=1000$.}
\label{fig:nu_k}
\end{figure}

The gradients $\nabla_{\theta}\mathcal{L}$ and $\nabla_{\theta_0}\mathcal{L}$ are computed with the \texttt{autodiff} algorithm, and the latter (the gradients with respect to the physical parameters) are scaled by $\gamma\nu(k)$ at each iteration. The scaling of the gradients is crucial for the convergence of the inverse PINN, as it ensures that the physical parameters are updated only when data is included in the loss function. The parameters are then updated with the Adam optimiser, with a sequence of learning rates that decrease at each iteration according to the epoch $k$. Namely starting from a learning rate $\alpha_0$ at the first epoch, the learning rate is updated as:
\begin{equation}
        \alpha_k = {\alpha_0} \beta^{\left\lfloor \frac{k}{100} \right\rfloor}
        \,,
\end{equation}
where $0.9<\beta<0.99$ is a constant factor. An algorithmic description of the above process is given in Algorithm \ref{alg:training}.

\begin{algorithm} 
\caption{Training Algorithm with Adaptive Weights and Gradient Updates.}
\label{alg:training} 

\begin{algorithmic}[1] 
\State $\mathrm{epoch}=0$
\Repeat
\State $\mathrm{epoch}=\mathrm{epoch}+1$
\If{$\mathrm{do\_parameter\_train}$ \textbf{and} $\mathrm{epoch}>K_0$}
\State compute $\nu(\mathrm{epoch})$ as in \textcolor{blue}{Eq.}\eqref{eq:nu_k} 
\EndIf
\State update data weights as in \textcolor{blue}{Eq.}\eqref{eq:omega}
\State compute gradients of loss function 
\State rescale gradients relative to $\theta_0$ by $\nu(\mathrm{epoch})$
\State apply gradients to all trainable parameters
\Until{convergence \textbf{or} $\mathrm{epoch}>\mathrm{epochs}$}
\end{algorithmic} 
\end{algorithm}

\section{Numerical results}
\label{sec:results}

In this section, we apply our PINN to different models arising from \cref{eq:adr} and \cref{eq:mi}, under several assumptions and conditions. All codes and data used in this work are published and freely available \citep{icardi_2024_12743929}. We report here a series of numerical experiments starting from a random initial guess for the parameters, and we show the convergence of the PINN to the correct values. The robustness of the approach with respect to the initial value of the parameters is shown in \cref{app:sensitivity} through two additional random initial values. Although it is outside the scope of this work to provide a quantitative sensitivity analysis, we show that the PINN is able to converge to the correct values for different initial guesses.

\subsection{Reference data}
\label{sec:refdata}
We use the \texttt{chebfun} package to generate reference data for the PINN training. The package is based on the Chebyshev polynomial approximation, and it is particularly suited for the solution of differential equations. We use the package to generate reference data for the PINN training, and we compare the results obtained with the PINN with the reference data generated by \texttt{chebfun}.

The numerical solution are computed at $N=100$ points in the spatial domain and $M=100$ points in the time domain, both uniformly sampled.  These are used both as training and collocation points.
Before performing the inverse PINN training, we have tested the PINN architecture for the direct problems, and we have verified that the PINN is able to accurately solve the direct problems for given parameter values. These have not been reported here for brevity.

\subsection{Pure diffusion}
\label{sec:testcase0}

The first testcase we consider is the pure diffusion problem, described by \cref{eq:heat}. We consider a one-dimensional spatial domain $\Omega=[0,1]$ and a time domain $[0,1]$. The initial condition is $u(x,0)=0$, and the boundary conditions are $u(0,t)=1$ and $\frac{\partial u}{\partial x}(1,t)=0$. This corresponds to a continuous source at the left boundary and a no-flux boundary condition at the right boundary.
The diffusion coefficient and its initial estimate are chosen randomly between $0.1$ and $10$. We use the PINN to solve the direct problem, and we consider the diffusion coefficient $D$ as a trainable parameter. We use the reference data generated by \texttt{chebfun} to train the PINN. 

We choose the following weights for the loss function: $\eta_{BC}=10$, $\eta_{IC}=10$, $\eta_{u}=1$. We use a total number of $L=9$ layers, i.e. eight hidden layers, and $m=20$ neurons for each hidden layer; moreover, we fix a total of $K=5000$ epochs, and we start updating the parameters and threshold epoch $K_0=1000$. The initial learning rate is set to $\alpha_0=0.01$ , the gradient scaling factor is set to $\gamma=0.2$ and learning rate reduction factor is set to $\beta=0.95$.

The computational time for the training of the PINN is approximately 800 seconds on an Apple Silicon M1 Pro processor with 10 cores and 16 GB of RAM.

\begin{figure}[htbp]
\centering
\includegraphics[width=0.49\textwidth]{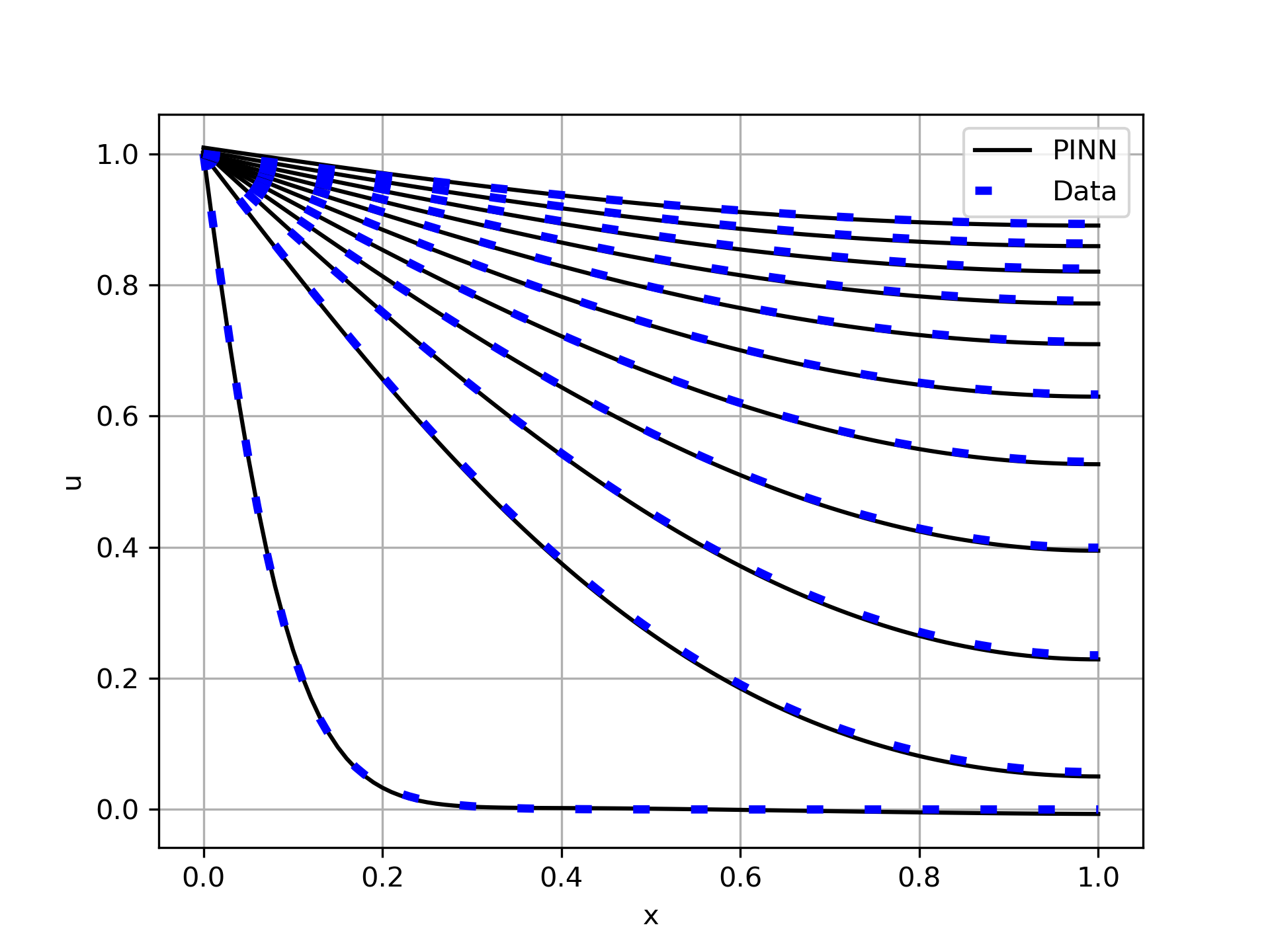}
\includegraphics[width=0.49\textwidth]{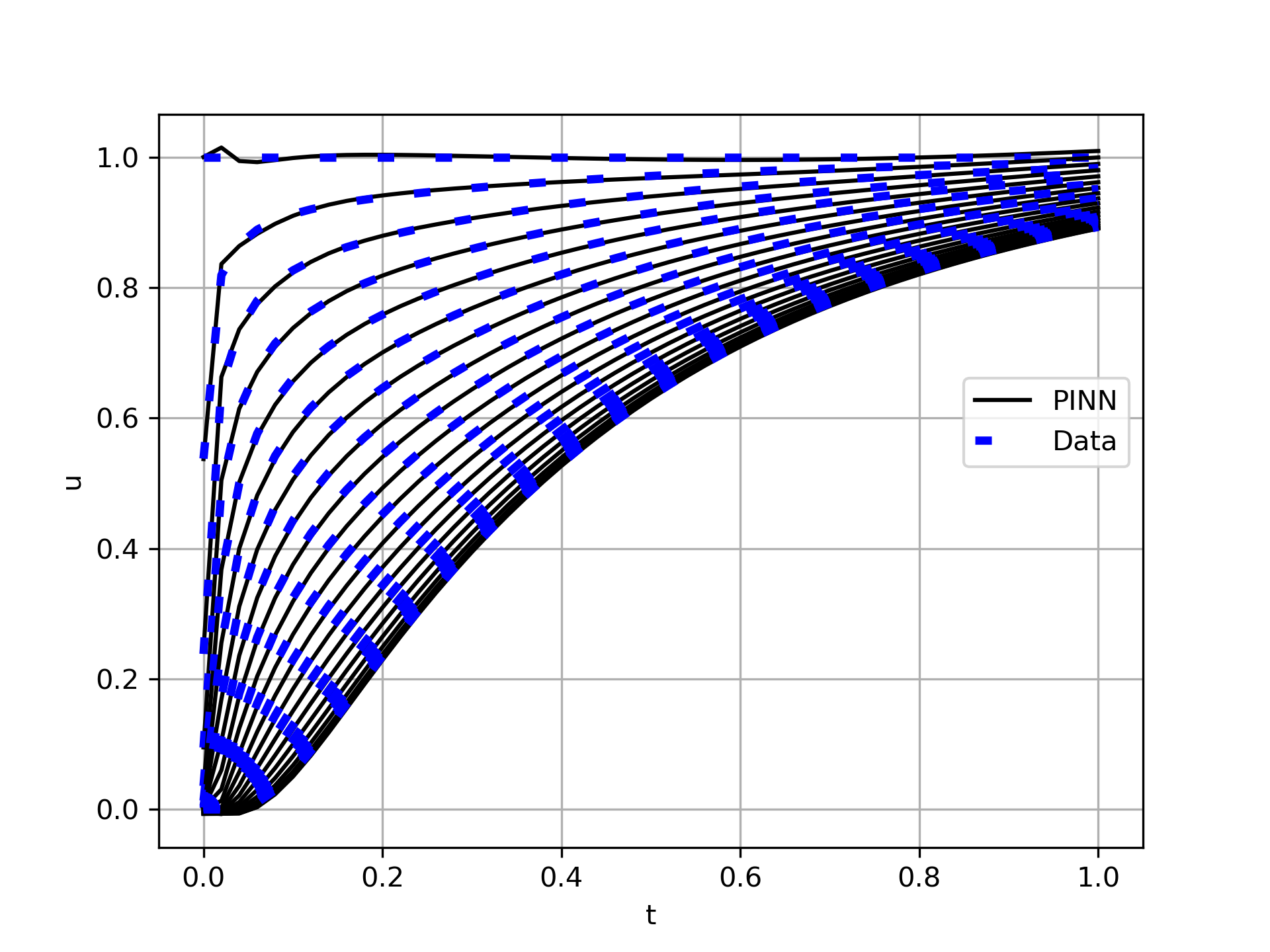}
\caption{Solution of the pure diffusion problem for the final parameter values. PINN approximation (continuous line) and reference data (dashed line). Concentration as a function of space for different times (left) and concentration as a function of time for different space locations (right).}
\label{fig:testcase0_sol}
\end{figure}

\begin{figure}[htbp]
\centering
\includegraphics[width=0.49\textwidth]{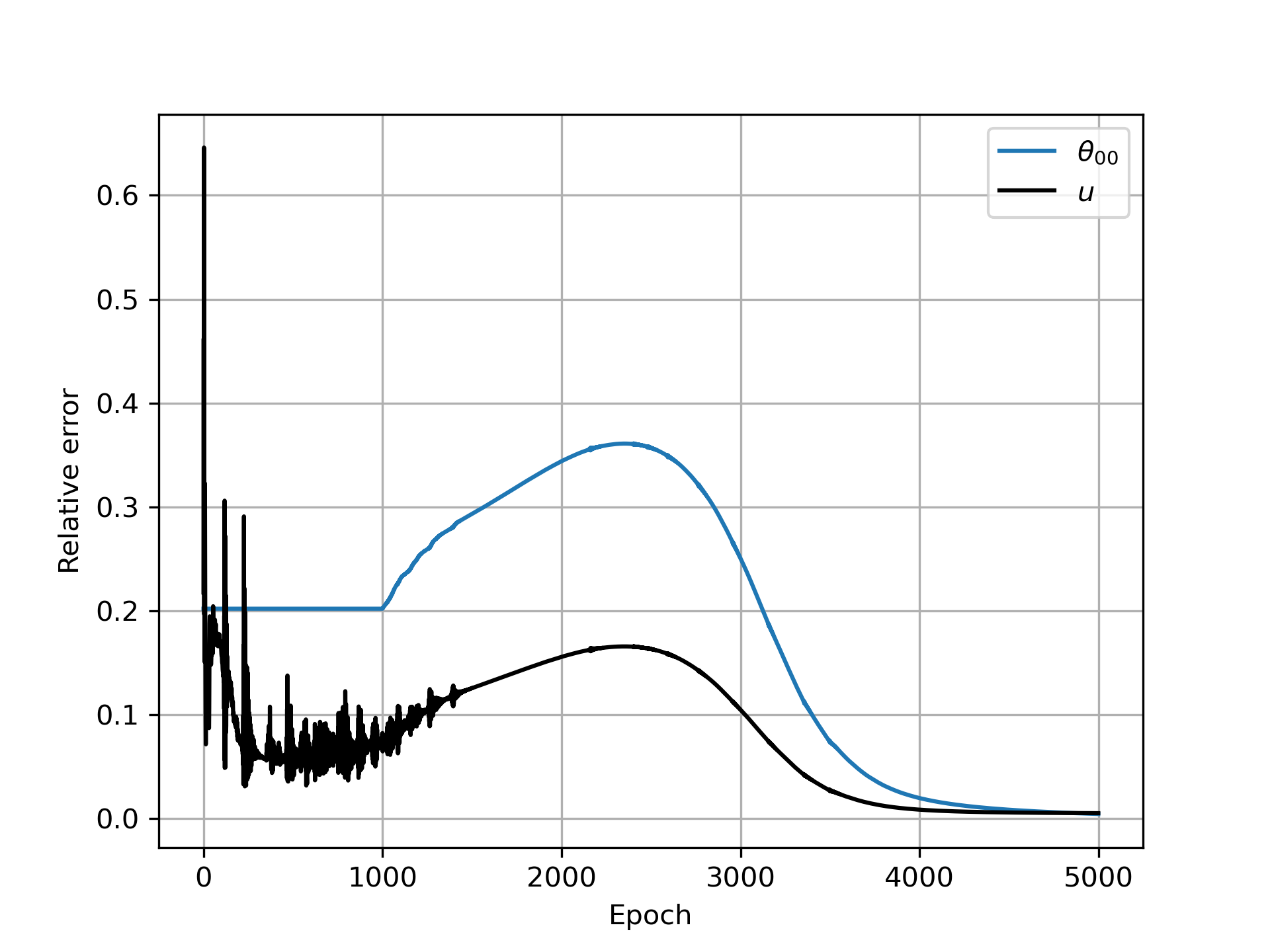}
\includegraphics[width=0.49\textwidth]{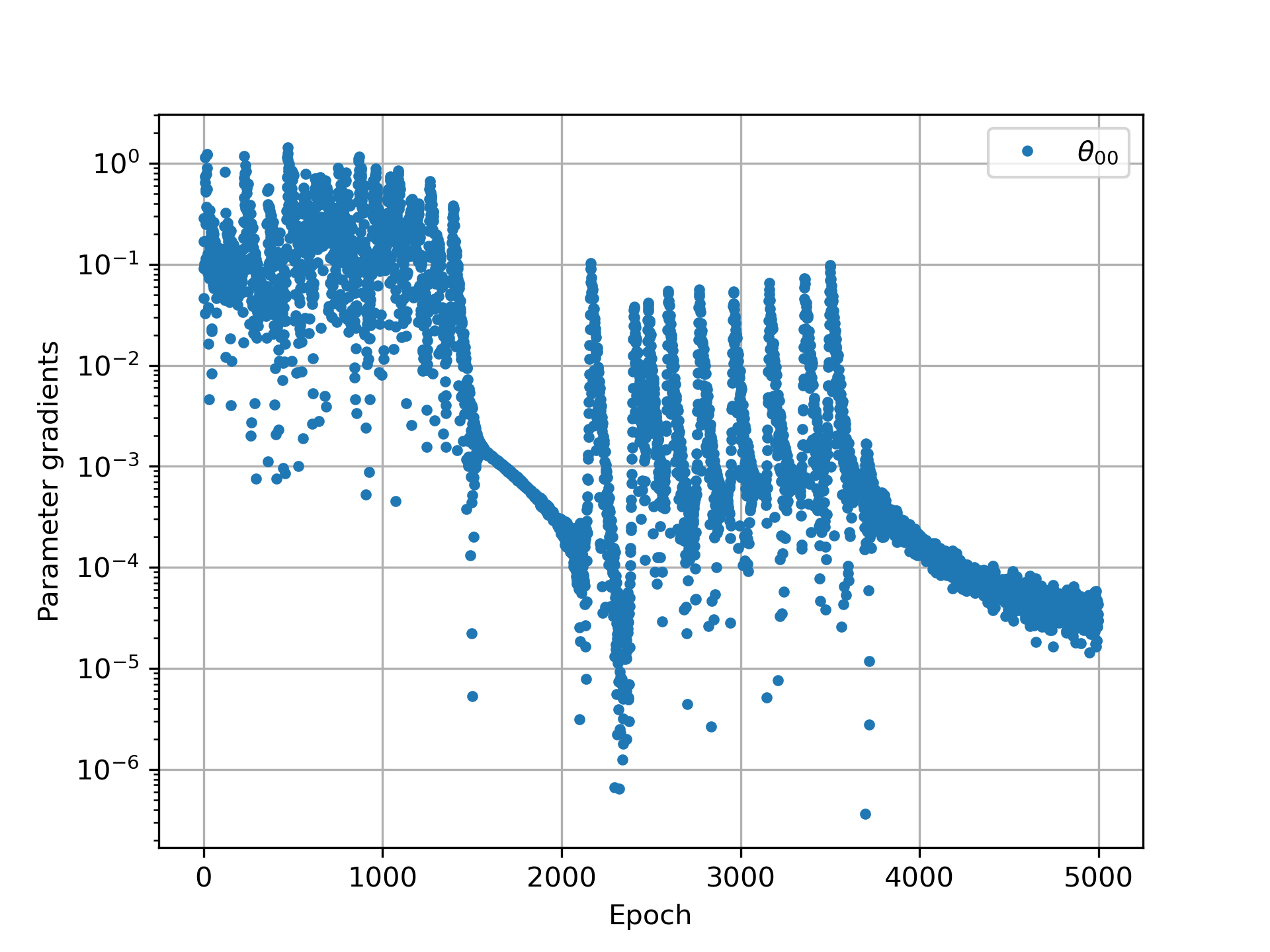}
\caption{Relative error for the diffusion coefficient and the solution $u$ during the training (left) and gradients of the diffusion coefficient (right).}
\label{fig:testcase0_param}
\end{figure}

\begin{figure}[htbp]
\centering
\includegraphics[width=0.49\textwidth]{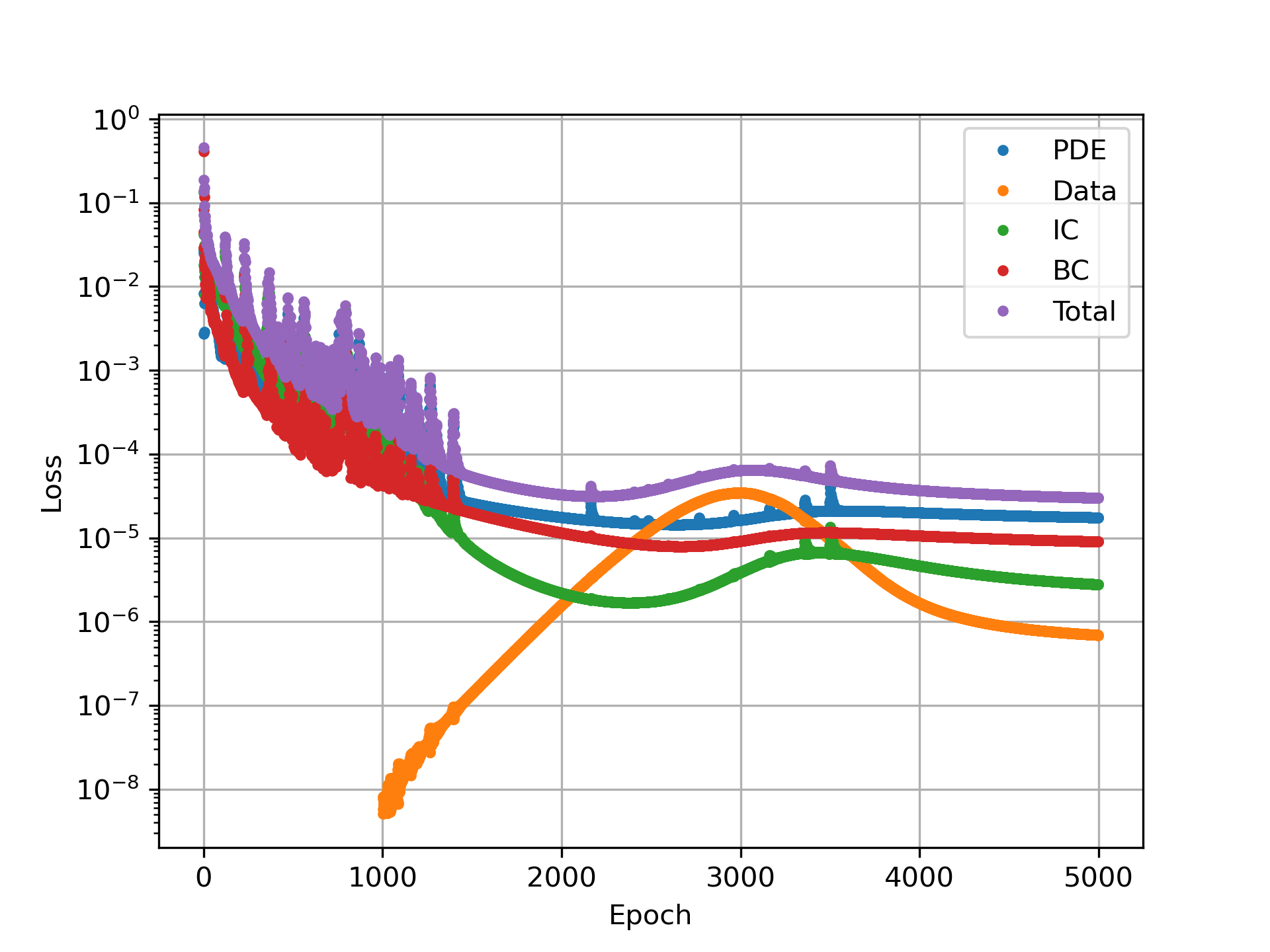}
\includegraphics[width=0.49\textwidth]{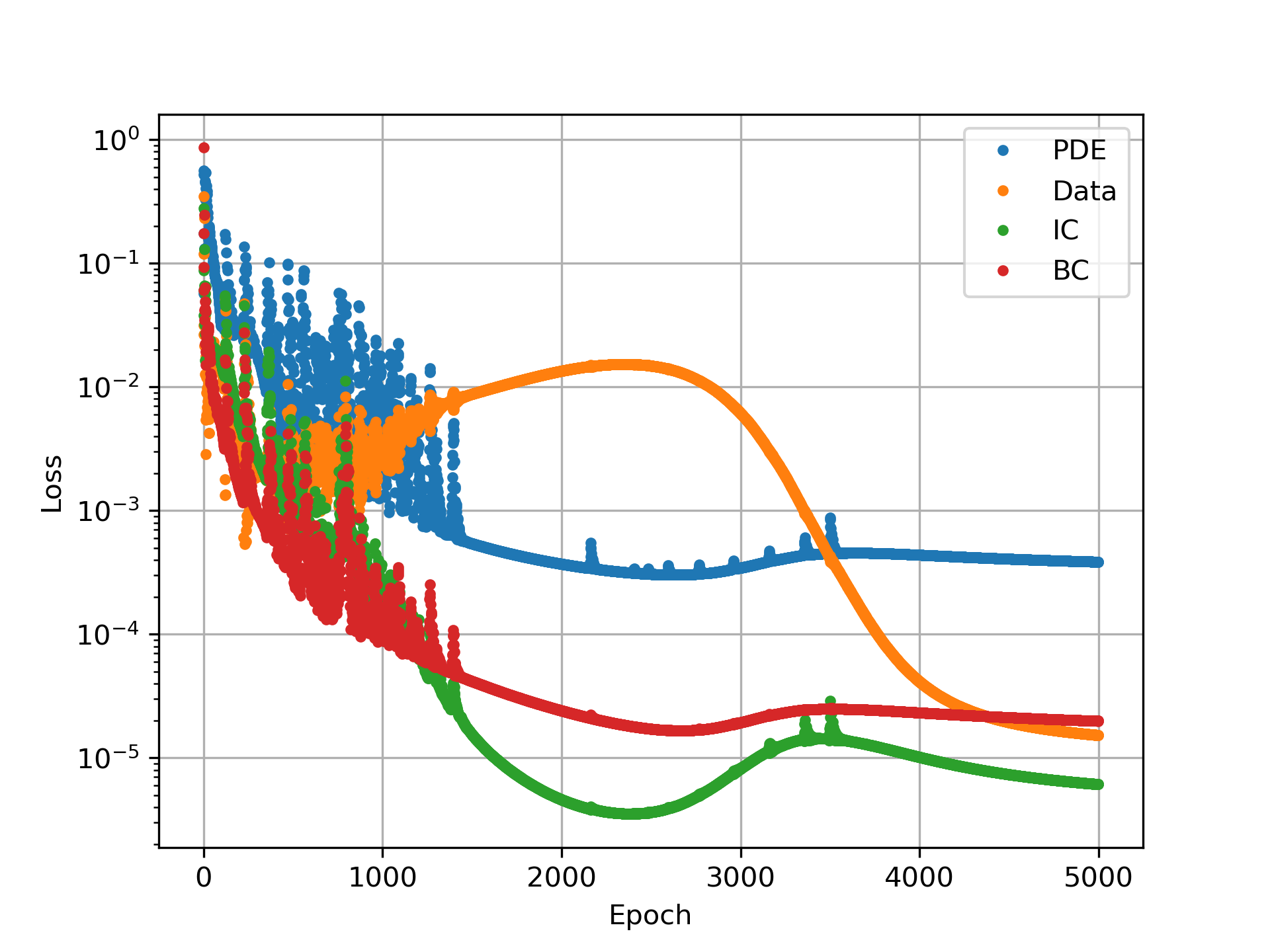}
\caption{Evolution of the weighted (left) and unweighted (right) loss functions during training of PINN.}
\label{fig:testcase0_loss}
\end{figure}

In \cref{fig:testcase0_sol}, we show the solution of the pure diffusion problem obtained with the PINN. The solution is compared with the reference data generated by \texttt{chebfun}. The solution is shown as a function of space for different times (left) and as a function of time for different space locations (right). The solution obtained with the PINN is in good agreement with the reference data, showing that the PINN is able to accurately solve the direct problem.

In \cref{fig:testcase0_param}, we show the evolution of the diffusion coefficient during the training of the PINN. The diffusion coefficient is shown as a function of the training iteration (epoch). The diffusion coefficient is updated during the training of the PINN, and it converges to the correct value. The gradients of the diffusion coefficient are also shown as a function of the training iteration (epoch).

In \cref{fig:testcase0_loss}, we show the evolution of the weighted and unweighted loss functions during the training of the PINN. The weighted loss function is shown as a function of the training iteration (epoch), and it is updated at each epoch to ensure that the different components of the loss function are balanced. The unweighted loss function is also shown as a function of the training iteration (epoch), to better highlight the effect of the weighting factors.


\subsection{Advection Diffusion}
\label{sec:testcase5}
This testcase considers the advection-diffusion problem described by \cref{eq:adr}. We consider a one-dimensional spatial domain $\Omega=[0,1]$ and a time domain $[0,1]$. The initial condition is $u(x,0)=0$, and the boundary conditions are $u(0,t)=H(0.01-t)$ and $\frac{\partial u}{\partial x}(1,t)=0$. This corresponds to a finite impulse at the left boundary and a no-flux boundary condition at the right boundary. The advection velocity and the dispersion coefficient are chosen randomly between 0.1 and 10, and the reaction term is here set to 0. We use the PINN to solve the direct problem, and we consider the advection velocity $V$ and the dispersion coefficient $D$ as trainable parameters. We use the reference data generated by \texttt{chebfun} to train the PINN.

We choose the following weights for the loss function: $\eta_{BC}=10$, $\eta_{IC}=10$, $\eta_{u}=2$. We use a total number of $L=9$ layers, i.e. eight hidden layers, and $m=20$ neurons for each hidden layer; moreover, we fix a total of $K=10000$ epochs, and we start updating the parameters and threshold epoch $K_0=1000$. The initial learning rate is set to $\alpha_0=0.01$ , the gradient scaling factor is set to $\gamma=0.2$ and learning rate reduction factor is set to $\beta=0.95$.

The computational time for the training of the PINN is approximately 1000 seconds on an Apple Silicon M1 Pro processor with 10 cores and 16 GB of RAM.

\begin{figure}[htbp]
\centering
\includegraphics[width=0.49\textwidth]{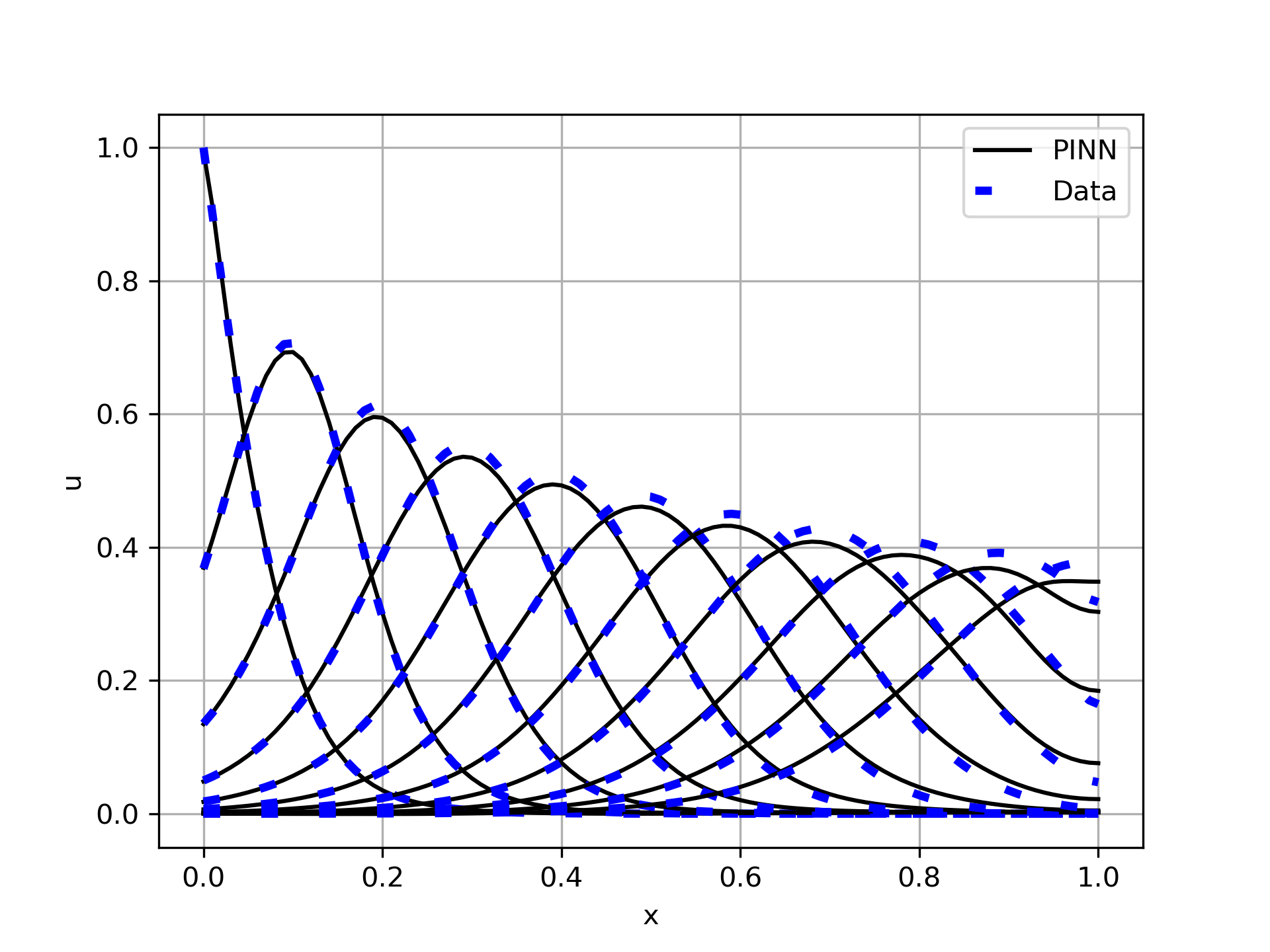}
\includegraphics[width=0.49\textwidth]{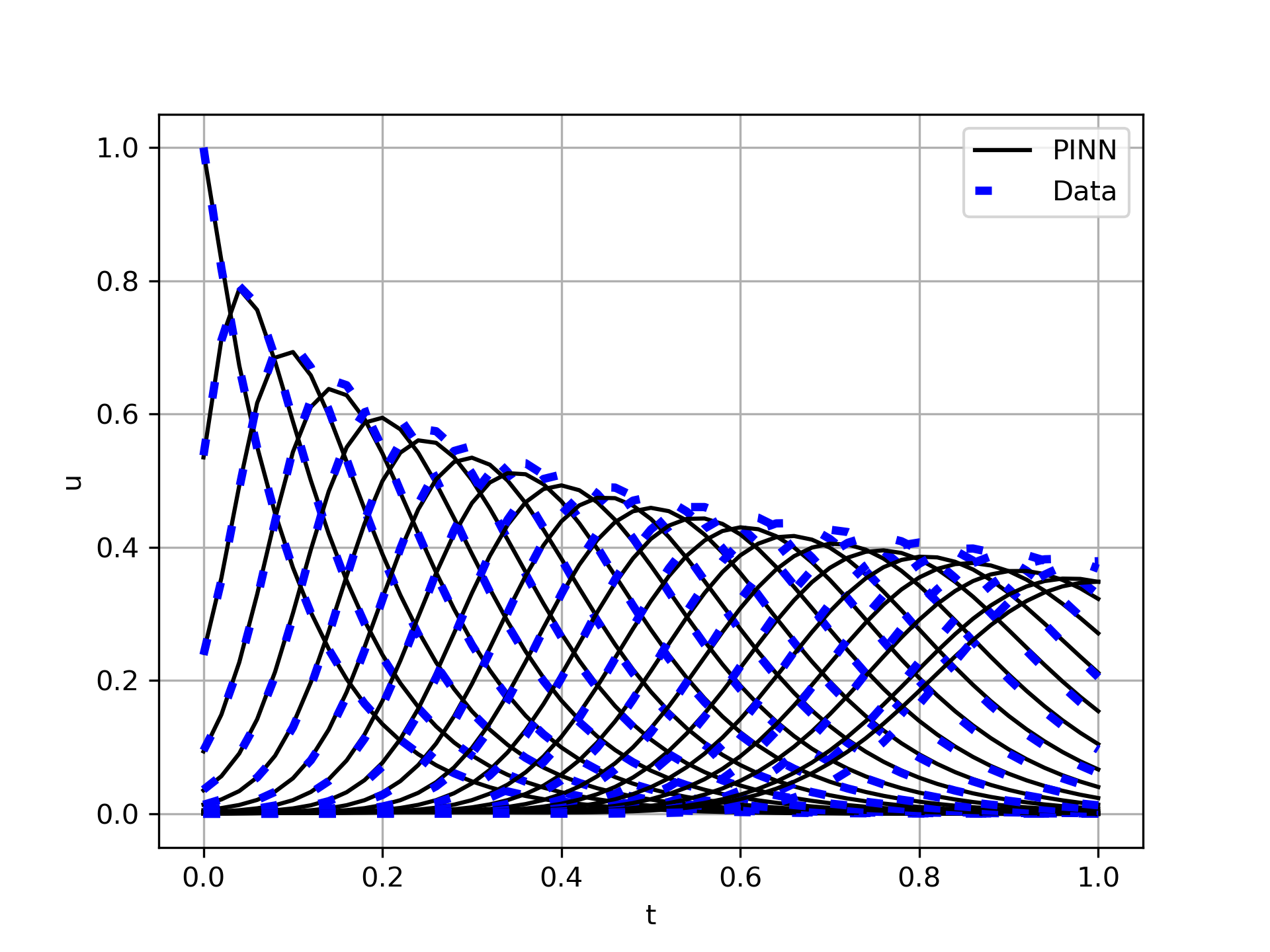}
\caption{Solution of the advection-diffusion problem for the final parameter values. PINN approximation (continuous line) and reference data (dashed line). Concentration as a function of space for different times (left) and concentration as a function of time for different space locations (right).}
\label{fig:testcase5_sol}
\end{figure}

\begin{figure}[htbp]
\centering
\includegraphics[width=0.49\textwidth]{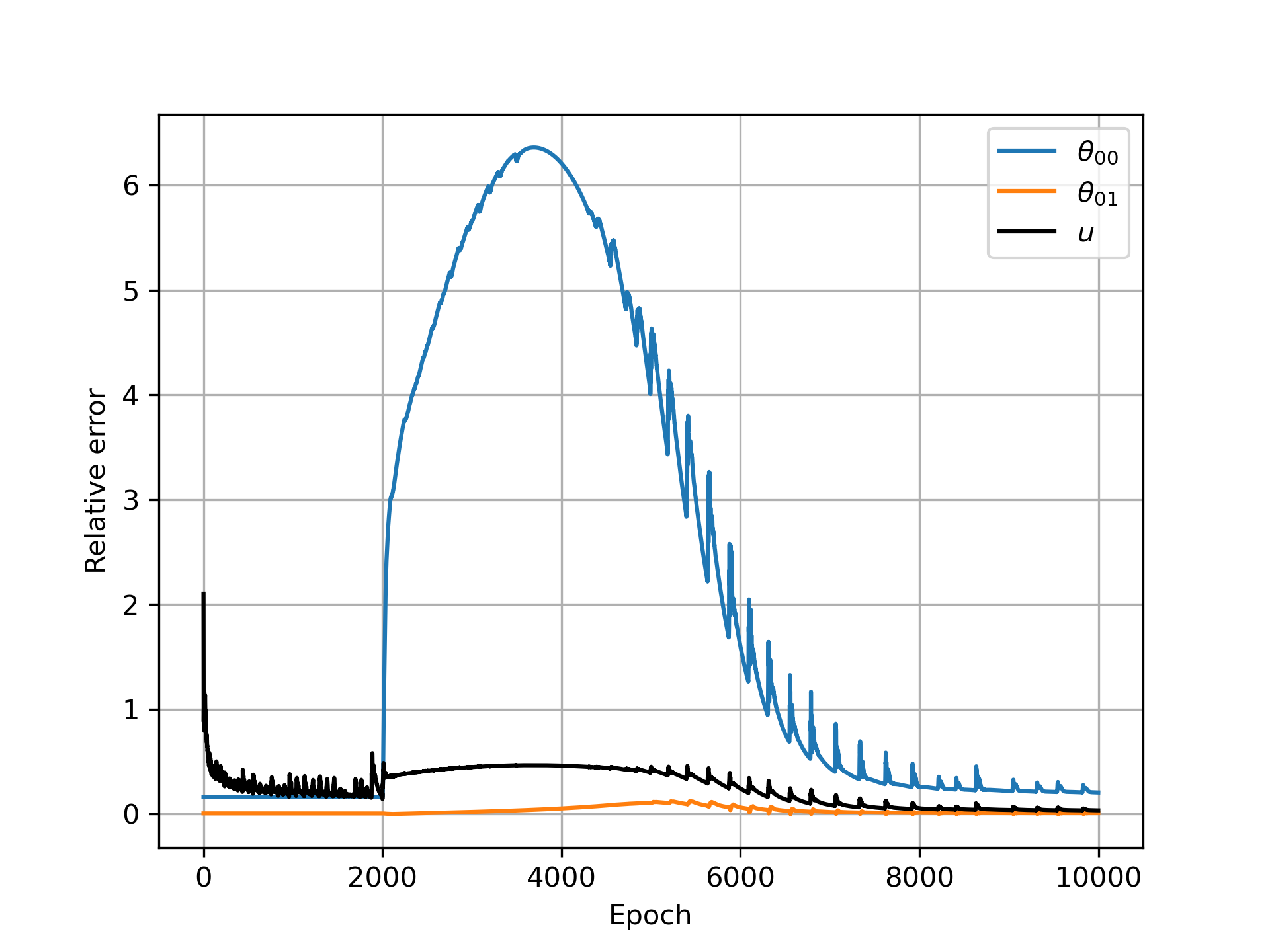}
\includegraphics[width=0.49\textwidth]{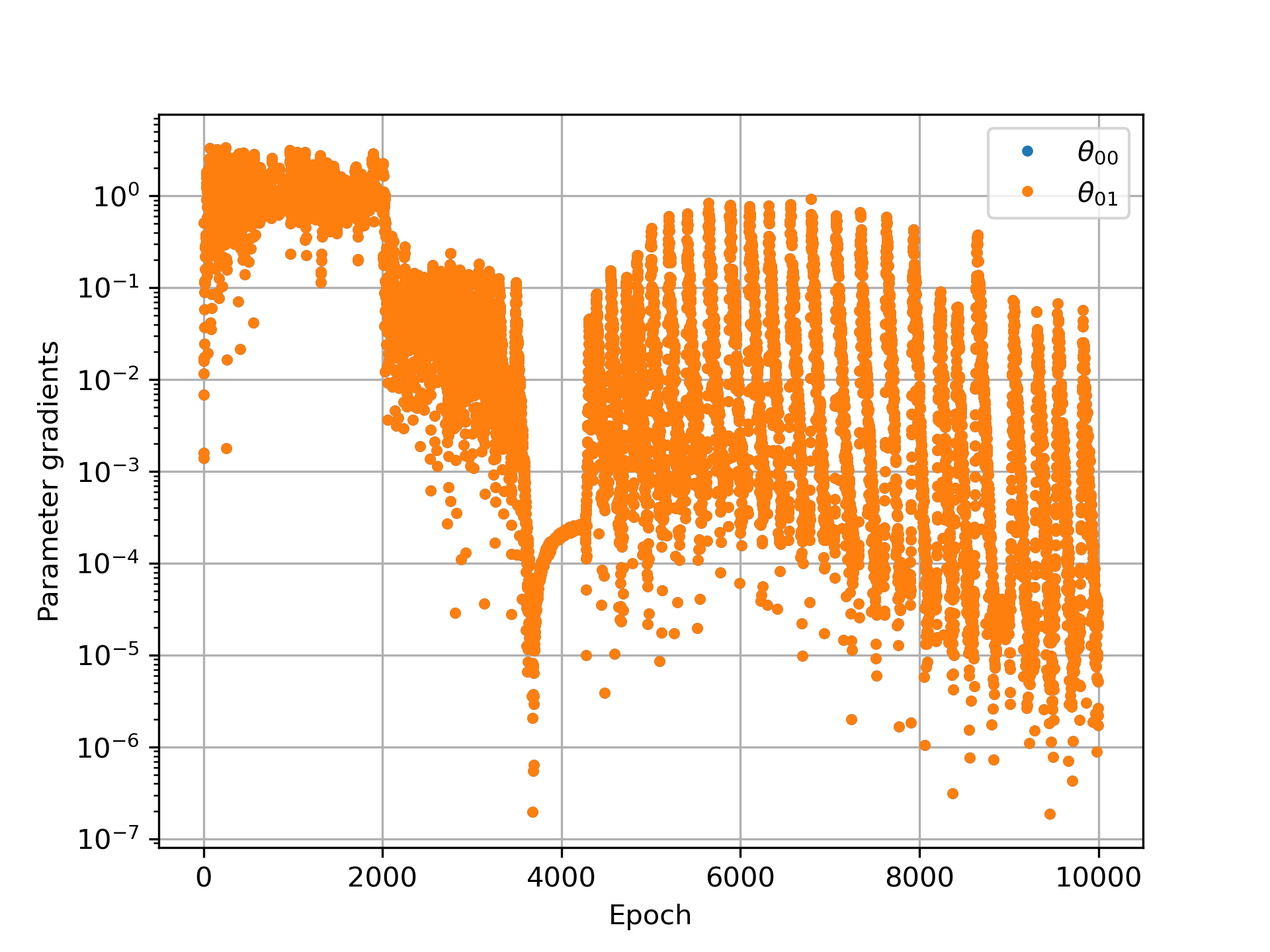}
\caption{Relative error for the advection velocity $\theta_{00}$, dispersion coefficient $\theta_{01}$ and the solution $u$ during the training (left) and gradients of the advection velocity and dispersion coefficient (right).}
\label{fig:testcase5_param}
\end{figure}

\begin{figure}[htbp]
\centering
\includegraphics[width=0.49\textwidth]{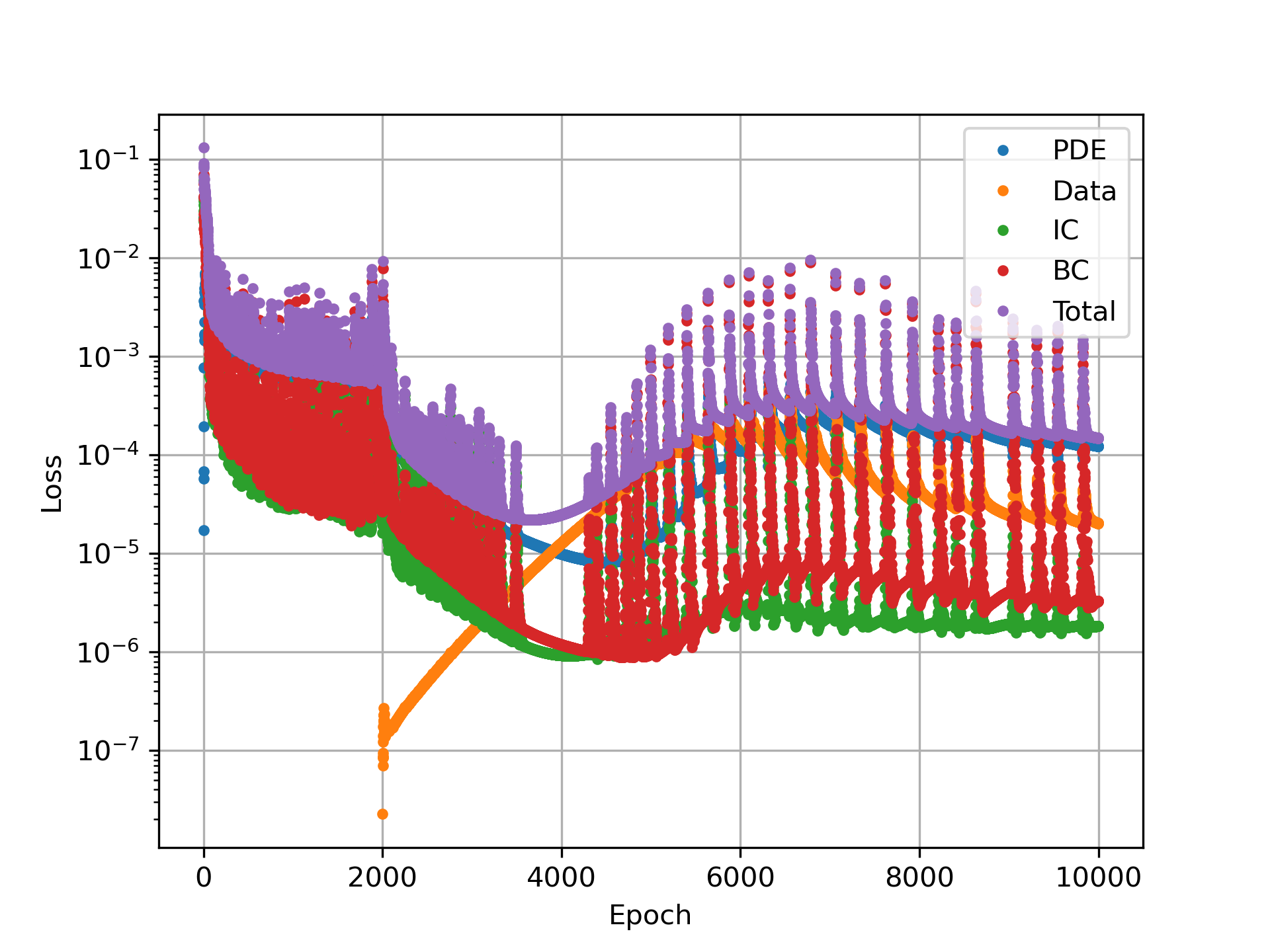}
\includegraphics[width=0.49\textwidth]{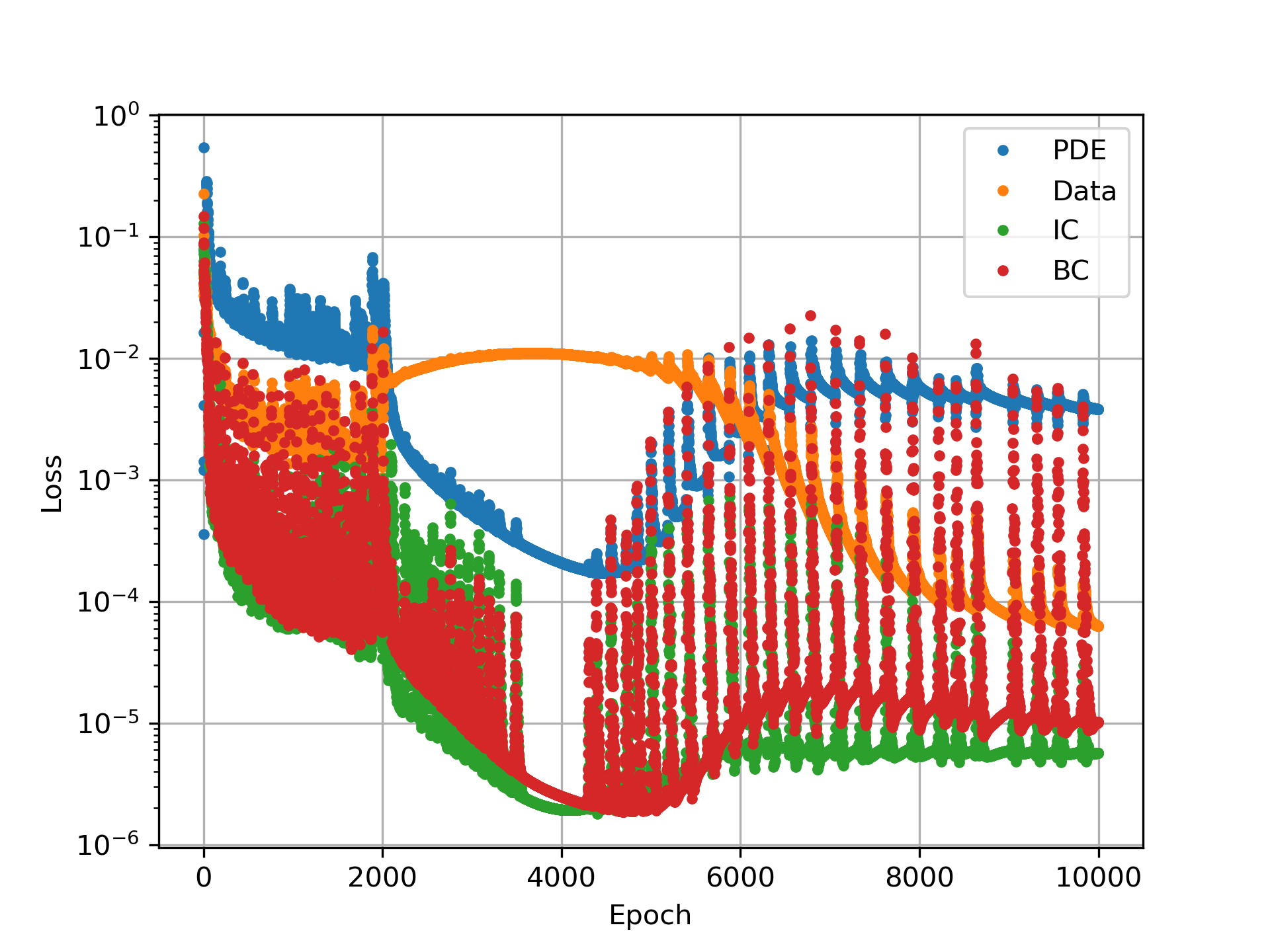}
\caption{Evolution of the weighted (left) and unweighted (right) loss functions during training of PINN.}
\label{fig:testcase5_loss}
\end{figure}

In \cref{fig:testcase5_sol}, we show the solution of the advection-diffusion problem obtained with the PINN. The solution is compared with the reference data generated by \texttt{chebfun}. The solution is shown as a function of space for different times (left) and as a function of time for different space locations (right). The solution obtained with the PINN is in good agreement with the reference data, showing that the PINN is able to accurately solve the direct problem.

In \cref{fig:testcase5_param}, we show the evolution of the advection velocity and dispersion coefficient during the training of the PINN. The advection velocity and dispersion coefficient are shown as a function of the training iteration (epoch). The advection velocity and dispersion coefficient are updated during the training of the PINN, and they converge to the correct values. The gradients of the advection velocity and dispersion coefficient are also shown as a function of the training iteration (epoch).

In \cref{fig:testcase5_loss}, we show the evolution of the weighted and unweighted loss functions during the training of the PINN. The weighted loss function is shown as a function of the training iteration (epoch), and it is updated at each epoch to ensure that the different components of the loss function are balanced. The unweighted loss function is also shown as a function of the training iteration (epoch), to better highlight the effect of the weighting factors.



\subsection{Advection-diffusion with mobile-immobile}
\label{sec:testcase1}

This testcase considers the advection-diffusion problem with mobile-immobile model described by \cref{eq:mi}. We consider a one-dimensional spatial domain $\Omega=[0,1]$ and a time domain $[0,1]$. The initial condition is $u(x,0)=v(x,0)=0$, and the boundary conditions are $u(0,t)=1$ and $\frac{\partial u}{\partial x}(1,t)=0$. This corresponds to a continuous injection at the left boundary and a no-flux boundary condition at the right boundary.
The effect of the immobile phase is to delay the transport of the solute, and the transfer coefficient $\lambda$ controls the transfer of solute between the mobile and immobile phases, resulting in long tails and non-Fickian transport behaviour \citep{municchi2021heterogeneous}.
The mobile and immobile porosities are fixed to $\beta_0=0.3$ and $\beta_1=0.1$, respectively, while dispersivity, effective velocity and transfer coefficient are the parameters to be estimated. The true values used for the training set are $D=0.1$, $V=1$, and $\lambda=10$. We use the PINN to solve the direct and inverse problem, and we consider the dispersion coefficient $D$, the effective velocity $V$, and the transfer coefficient $\lambda$ as trainable parameters. We use the reference data generated by \texttt{chebfun} to train the PINN. The initial values for the parameters are chosen randomly between 0.1 and 10.

We choose the following weights for the loss function: $\eta_{BC}=10$, $\eta_{IC}=10$, $\eta_{u}=1$. We use a total number of $L=11$ layers, i.e. ten hidden layers, and $m=25$ neurons for each hidden layer; moreover, we fix a total of $K=10000$ epochs, and we start updating the parameters and threshold epoch $K_0=2000$. The initial learning rate is set to $\alpha_0=0.01$ , the gradient scaling factor is set to $\gamma=0.1$ and learning rate reduction factor is set to $\beta=0.98$.

The computational time for the training of the PINN is approximately 2000 seconds on an Apple Silicon M1 Pro processor with 10 cores and 16 GB of RAM.

\begin{figure}[htbp]
\centering
\includegraphics[width=0.49\textwidth]{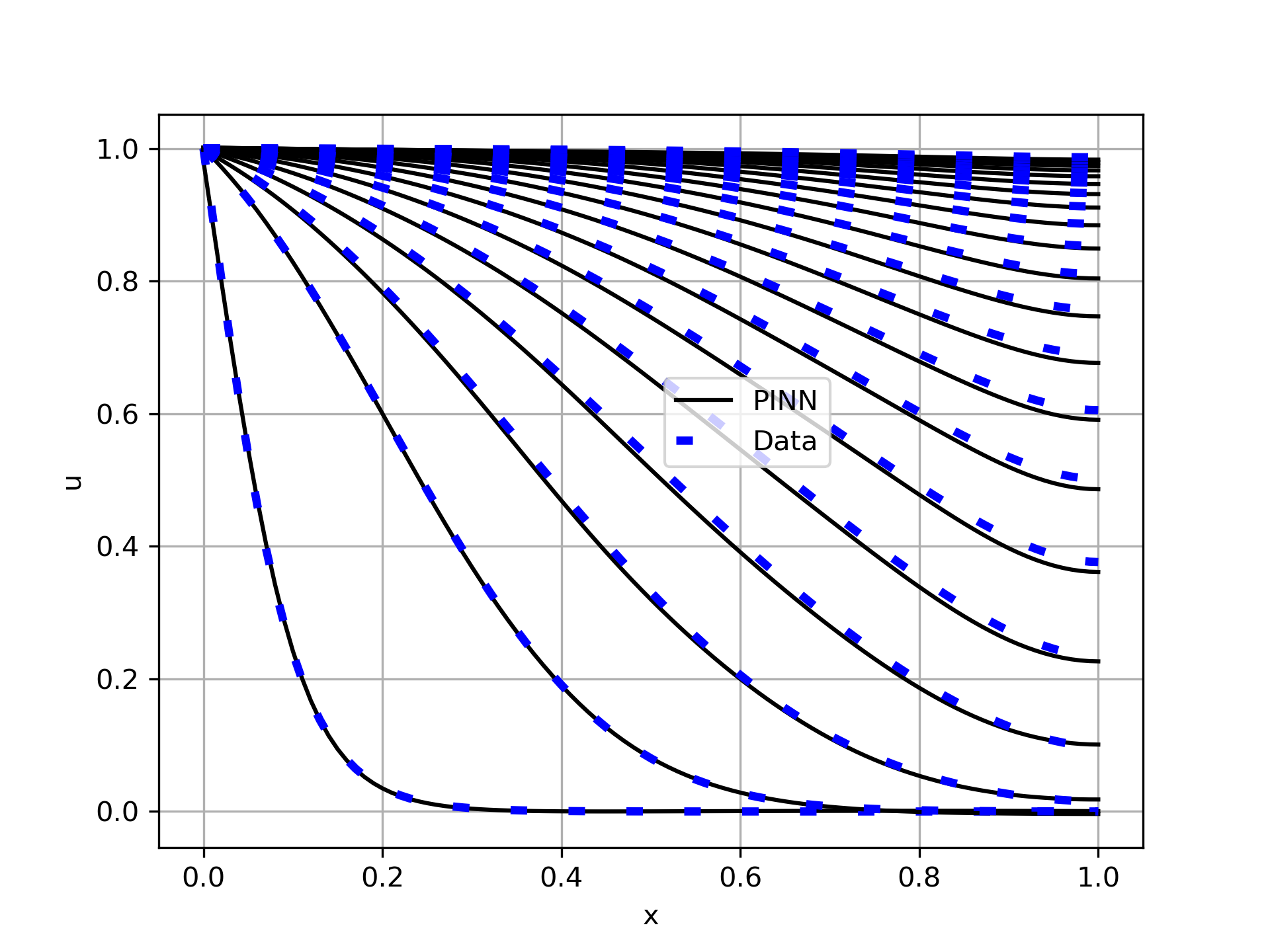}
\includegraphics[width=0.49\textwidth]{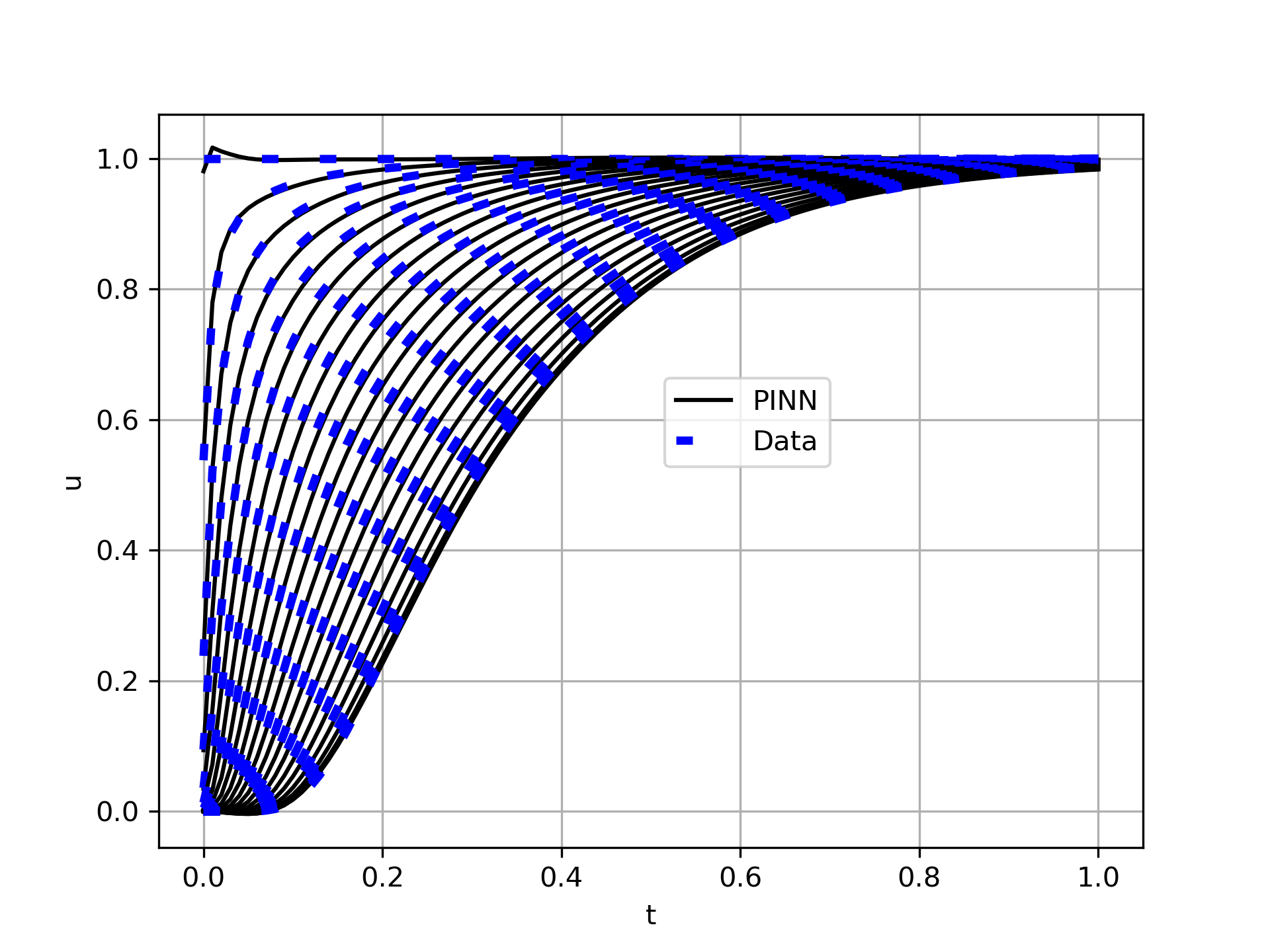}
\caption{Solution of the advection-diffusion with mobile-immobile model for the final parameter values. Mobile concentration. PINN approximation (continuous line) and reference data (dashed line). Concentration as a function of space for different times (left) and concentration as a function of time for different space locations (right).}
\label{fig:testcase1_sol}
\end{figure}

\begin{figure}[htbp]
        \centering
        \includegraphics[width=0.49\textwidth]{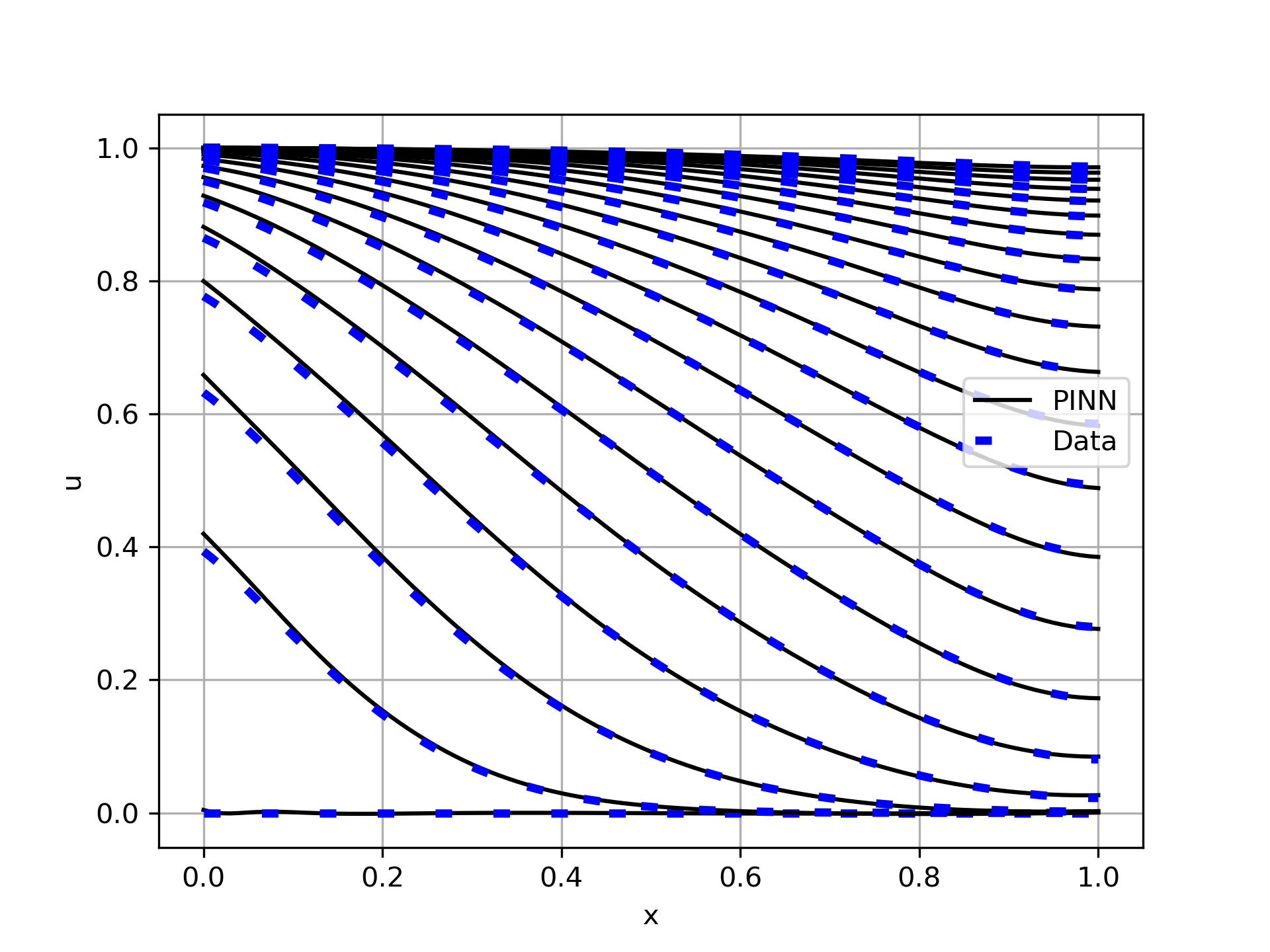}
        \includegraphics[width=0.49\textwidth]{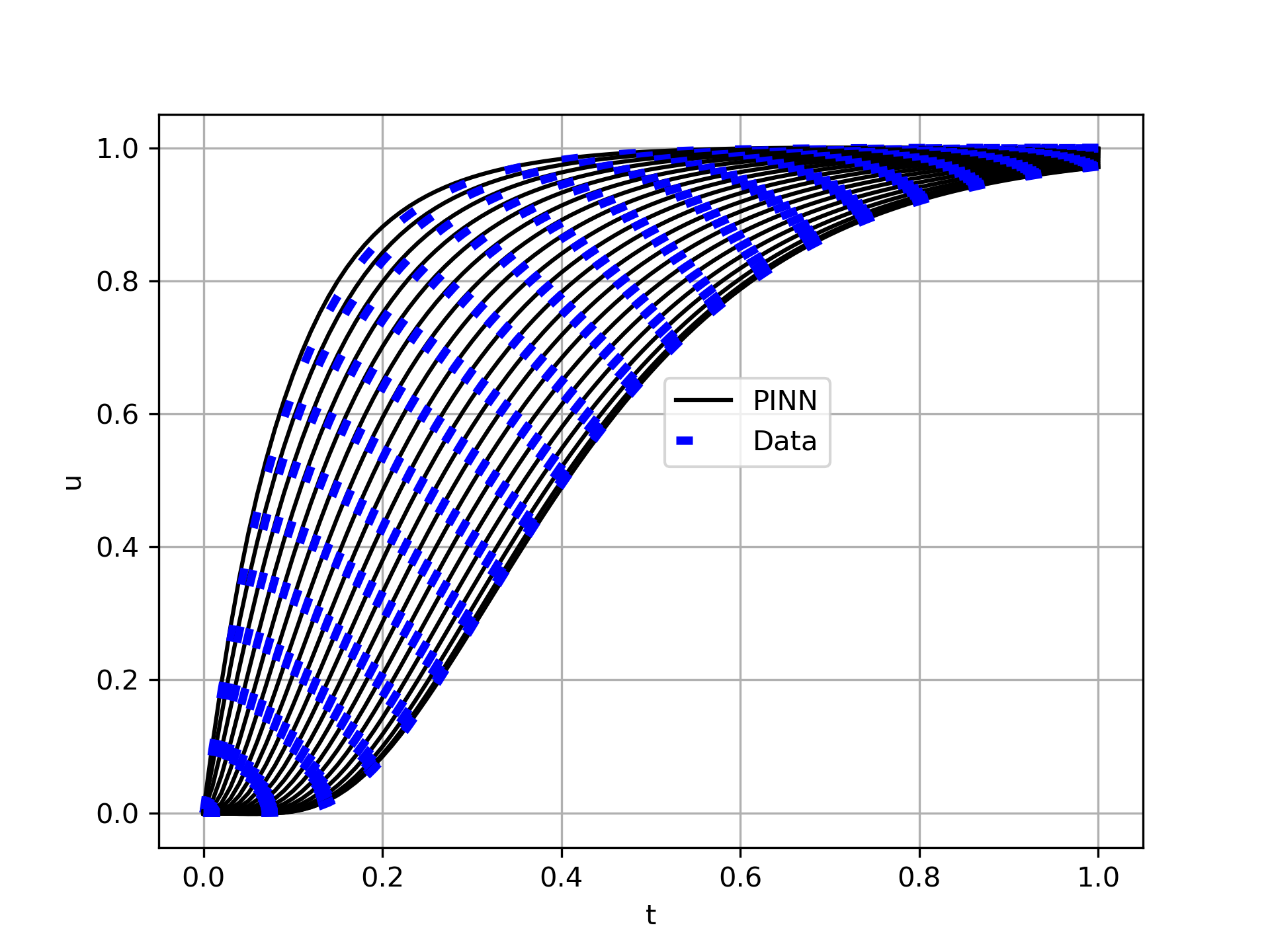}
        \caption{Solution of the advection-diffusion with mobile-immobile model for the final parameter values. Immobile concentration. PINN approximation (continuous line) and reference data (dashed line). Concentration as a function of space for different times (left) and concentration as a function of time for different space locations (right).}
        \label{fig:testcase1_sol1}
\end{figure}

\begin{figure}[htbp]
\centering
\includegraphics[width=0.49\textwidth]{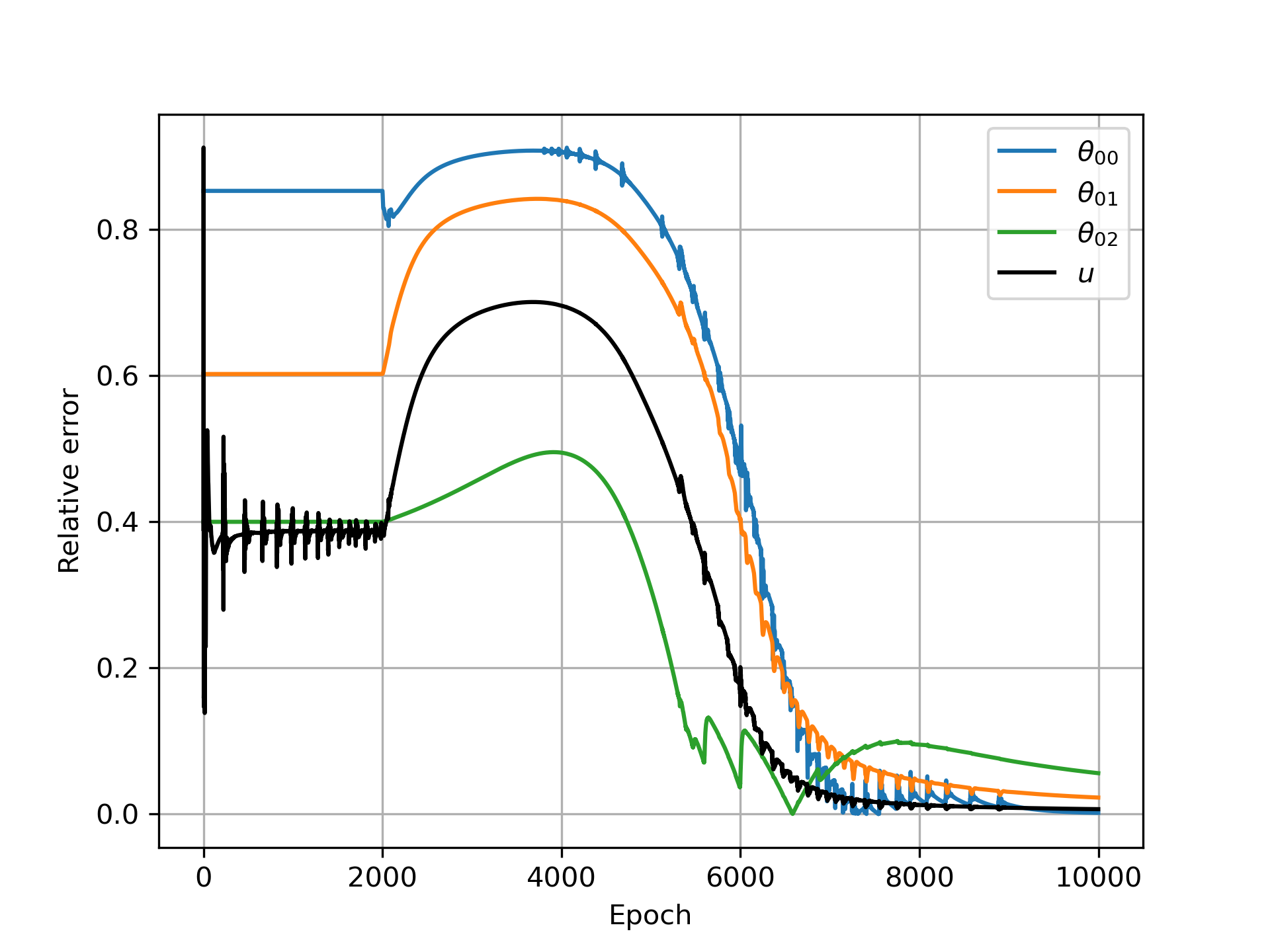}
\includegraphics[width=0.49\textwidth]{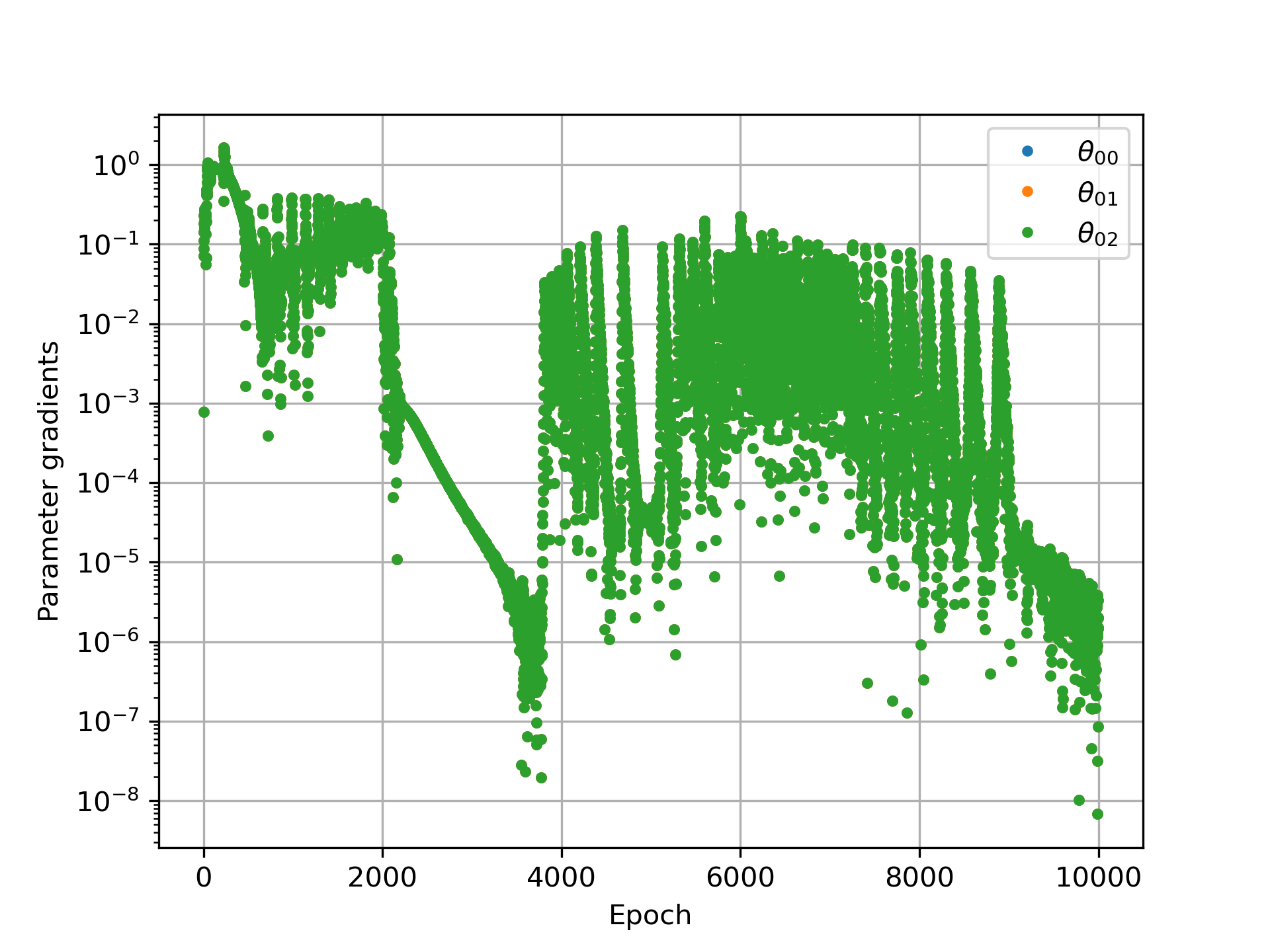}
\caption{Relative error for the dispersion coefficient $\theta_{00}$, advection velocity $\theta_{01}$, and transfer coefficient $\theta_{02}$ and the solution $u$ during the training (left) and gradients of the dispersion coefficient, effective velocity, and transfer coefficient (right).}
\label{fig:testcase1_param}
\end{figure}

\begin{figure}[htbp]
\centering
\includegraphics[width=0.49\textwidth]{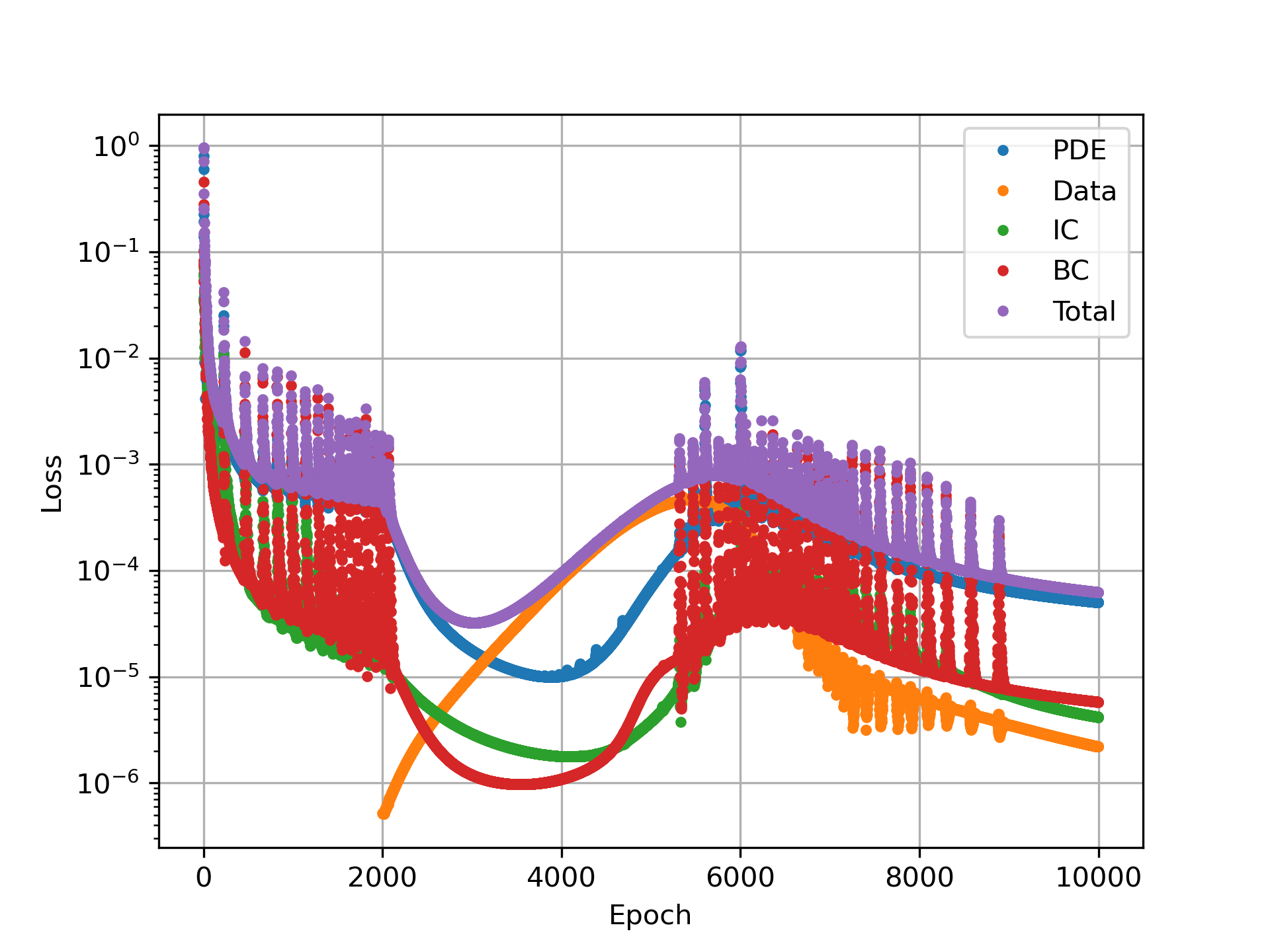}
\includegraphics[width=0.49\textwidth]{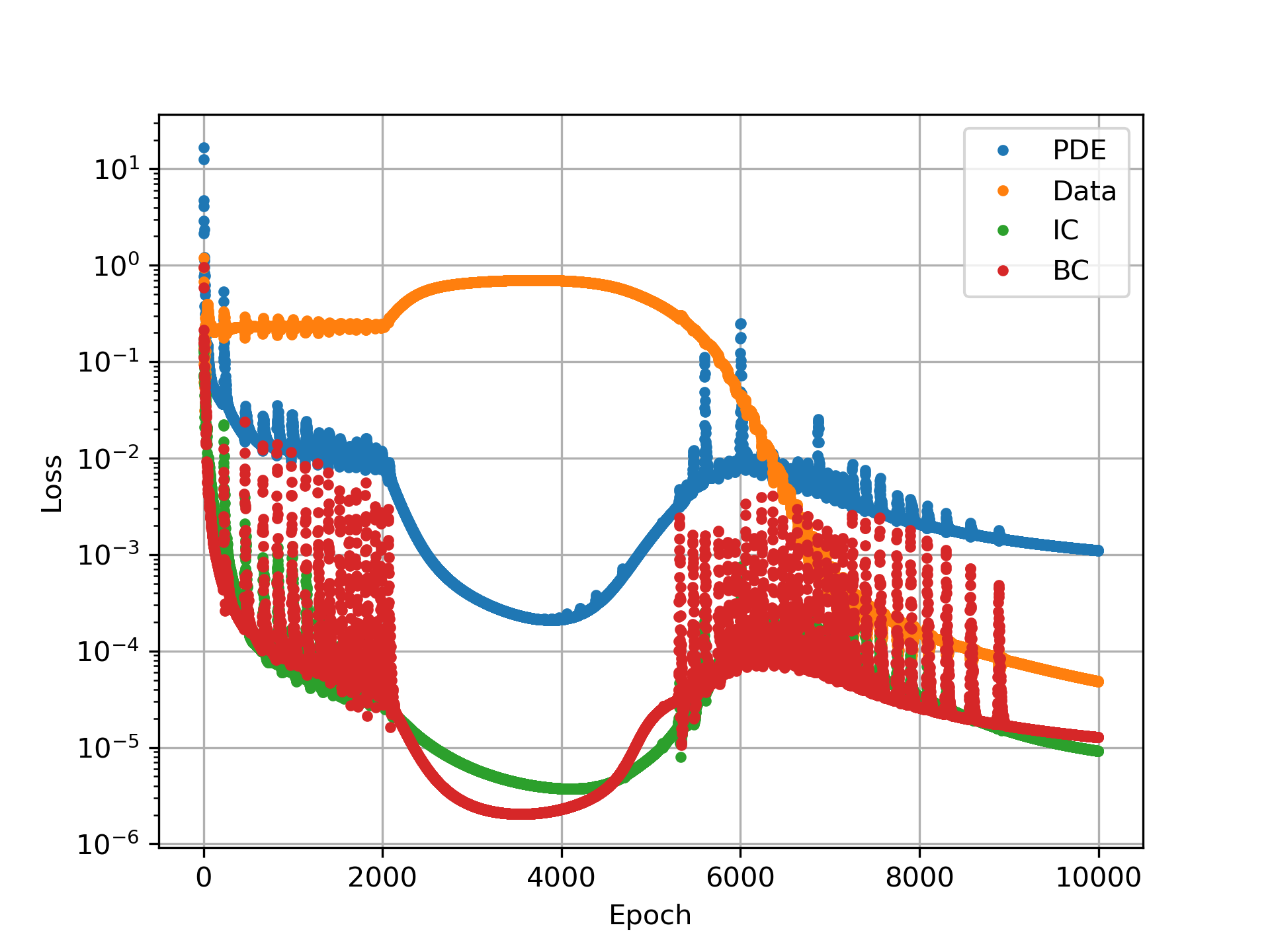}
\caption{Evolution of the weighted (left) and unweighted (right) loss functions during training of PINN.}
\label{fig:testcase1_loss}
\end{figure}

In \cref{fig:testcase1_sol,fig:testcase1_sol1} we show the solution of the advection-diffusion with mobile-immobile model obtained with the PINN, for the mobile and immobile phase, respectively. The solution is compared with the reference data generated by \texttt{chebfun}. The solution is shown as a function of space for different times (left) and as a function of time for different space locations (right). The solution obtained with the PINN is in good agreement with the reference data, showing that the PINN is able to accurately solve the direct problem. 

In \cref{fig:testcase1_param}, we show the evolution of the dispersion coefficient, effective velocity, and transfer coefficient during the training of the PINN. The relative error in the dispersion coefficient, effective velocity, and transfer coefficient is shown as a function of the training iteration (epoch). The parameters are updated during the training of the PINN, and they converge to the correct values. The gradients of the dispersion coefficient, effective velocity, and transfer coefficient are also shown as a function of the training iteration (epoch). 

In \cref{fig:testcase1_loss}, we show the evolution of the weighted and unweighted loss functions during the training of the PINN. The weighted loss function is shown as a function of the training iteration (epoch), and it is updated at each epoch to ensure that the different components of the loss function are balanced. The unweighted loss function is also shown as a function of the training iteration (epoch), to better highlight the effect of the weighting factors.

\section{Conclusions}
\label{sec:conclusions}

In this paper, we have proposed an adaptive inverse PINN architecture for solving transport problems in porous materials. These include a diffusion, advection-diffusion and mobile-immobile formulations. We propose a robust PINN architecture and training algorithm that can reproduce well the forward problem and the inverse problem for up to three parameters. Ongoing work include the extension to two-dimensional problems, larger number of parameters and non-parametric functional form of the parameters to include heterogeneities and non-linear dependencies: moreover, a comparison with other methods -data assimilation, Bayesian techniques- for handling inverse models will be carried out in future works; in a real world application framework, it could be of interest estimating some hydraulic parameters from real data, for instance in a highly nonlinear model such as Richards' equation, starting from soil water content measures over time (e.g., \cite{DeCarlo2018}).  The main novelty of the proposed approach is the adaptive scaling of the loss function components and gradients of the trainable parameters. This adaptive scaling is crucial for the convergence of the inverse problem, as it ensures that the different components of the loss function are balanced and that the gradients of the trainable parameters are scaled appropriately. We have demonstrated the effectiveness of the proposed approach through a series of numerical experiments, showing that the adaptive inverse PINN architecture is scalable, robust, and efficient for solving a wide range of transport models.

\section*{Acknowledgements}
MB and FVD are part of INdAM research group GNCS. \\MI gratefully acknowledges the support of Short Term Mobility program at  CNR,
funded by the Italian Ministry of Universities and Research (MUR) and  by the European Union through Next Generation EU, M4C2. \\
MB and FVD acknowledge the support of PRIN2022PNRR n. P2022M7JZW \emph{SAFER MESH - Sustainable mAnagement oF watEr Resources ModEls and numerical MetHods} research grant, funded by the Italian Ministry of Universities and Research (MUR) and  by the European Union through Next Generation EU, M4C2, CUP H53D23008930001. \\
MB thanks Mrs Domenica Livorti from CNR-IRSA for supporting the project activities.

\bibliographystyle{unsrtnat} 
\bibliography{nonFickianPINN}

\appendix

\section{Sensitivity with respect to the initial parameters}
\label{app:sensitivity}

\begin{figure}[htbp]
\centering
\includegraphics[width=0.49\textwidth]{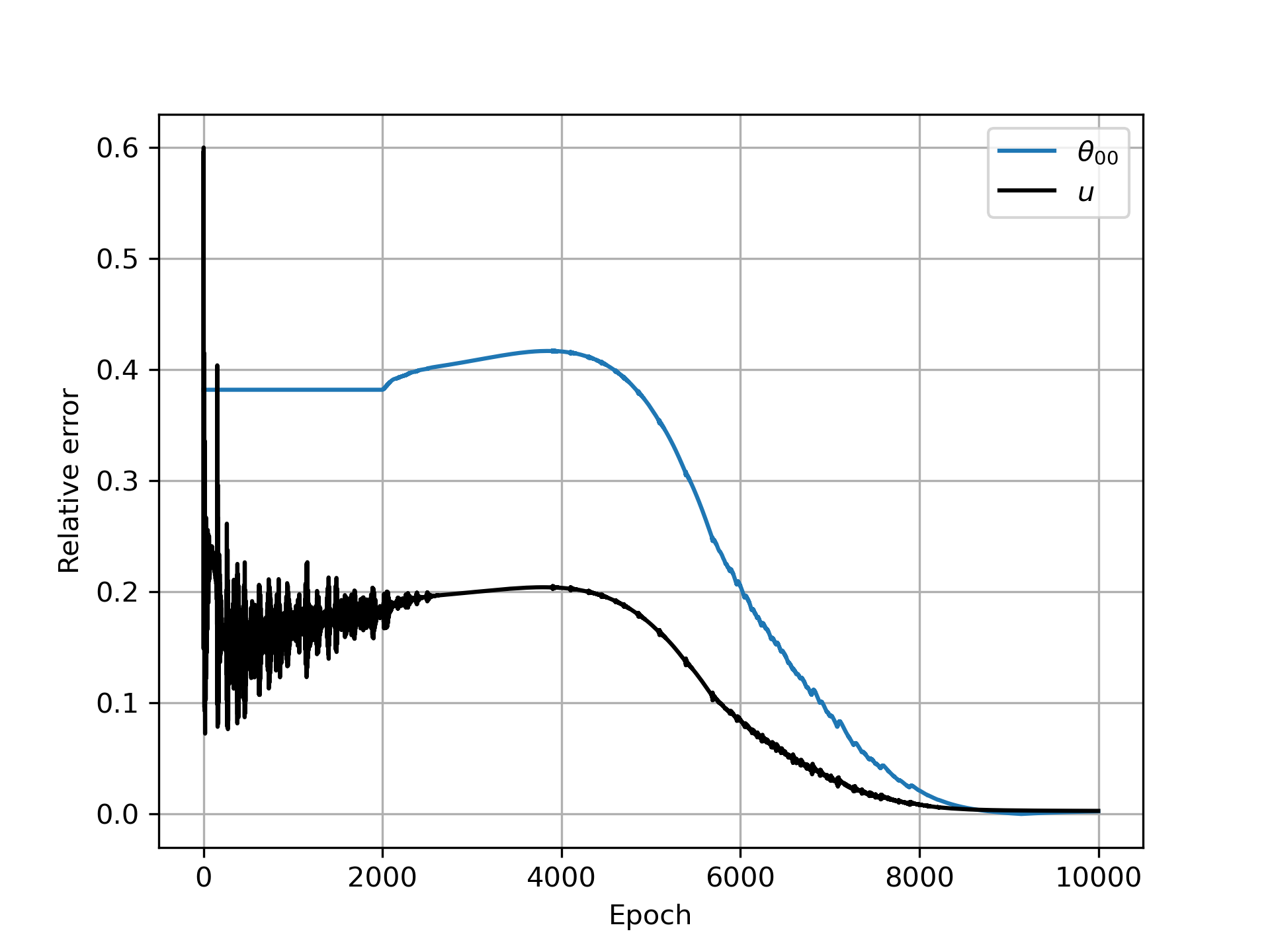}
\includegraphics[width=0.49\textwidth]{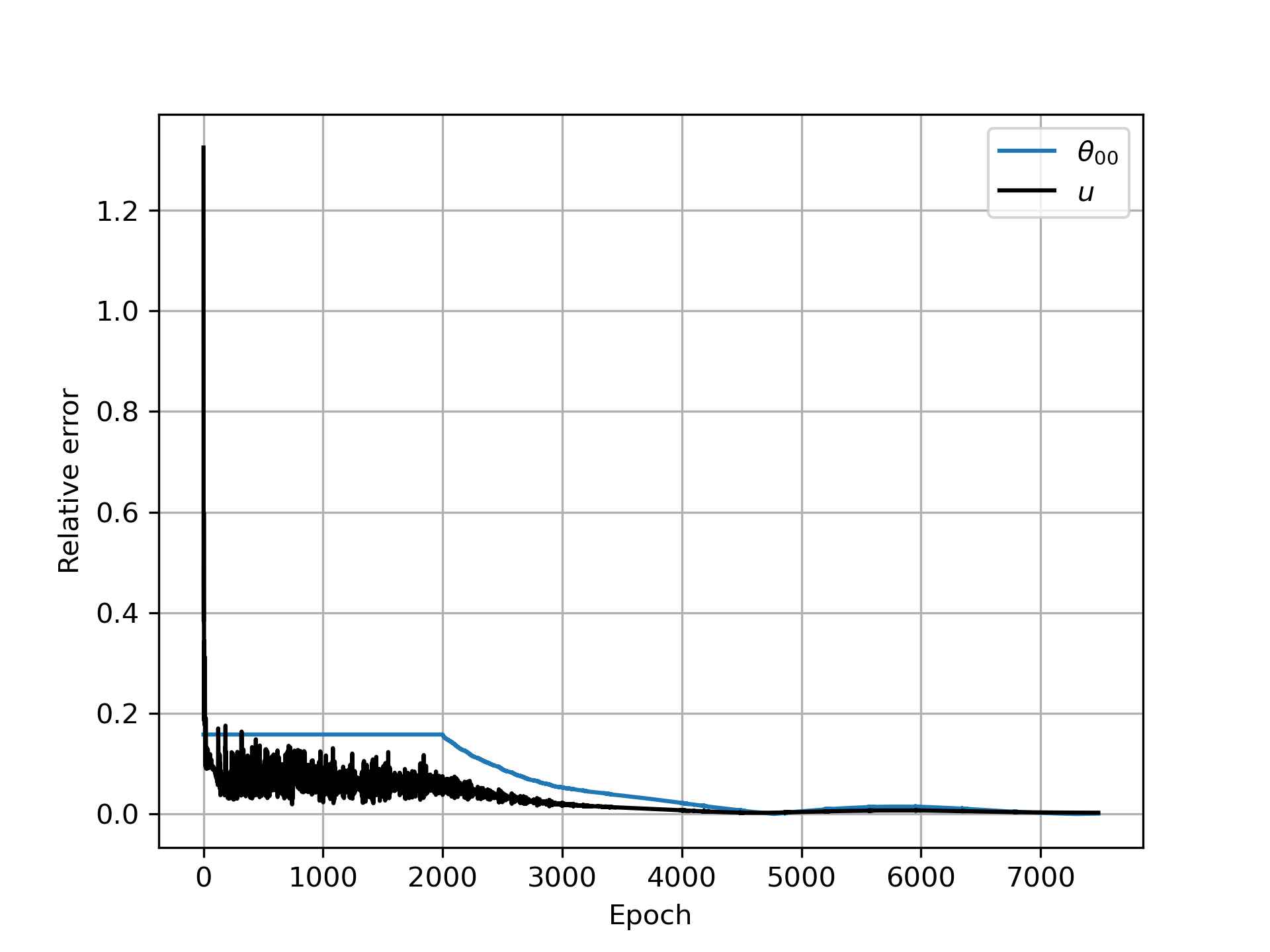}
\caption{Pure diffusion testcase. Relative error for the diffusion coefficient and the solution $u$ during the training for two additional random initial value of the parameters.}
\label{fig:appendix_param0}
\end{figure}

\begin{figure}[htbp]
\centering
\includegraphics[width=0.49\textwidth]{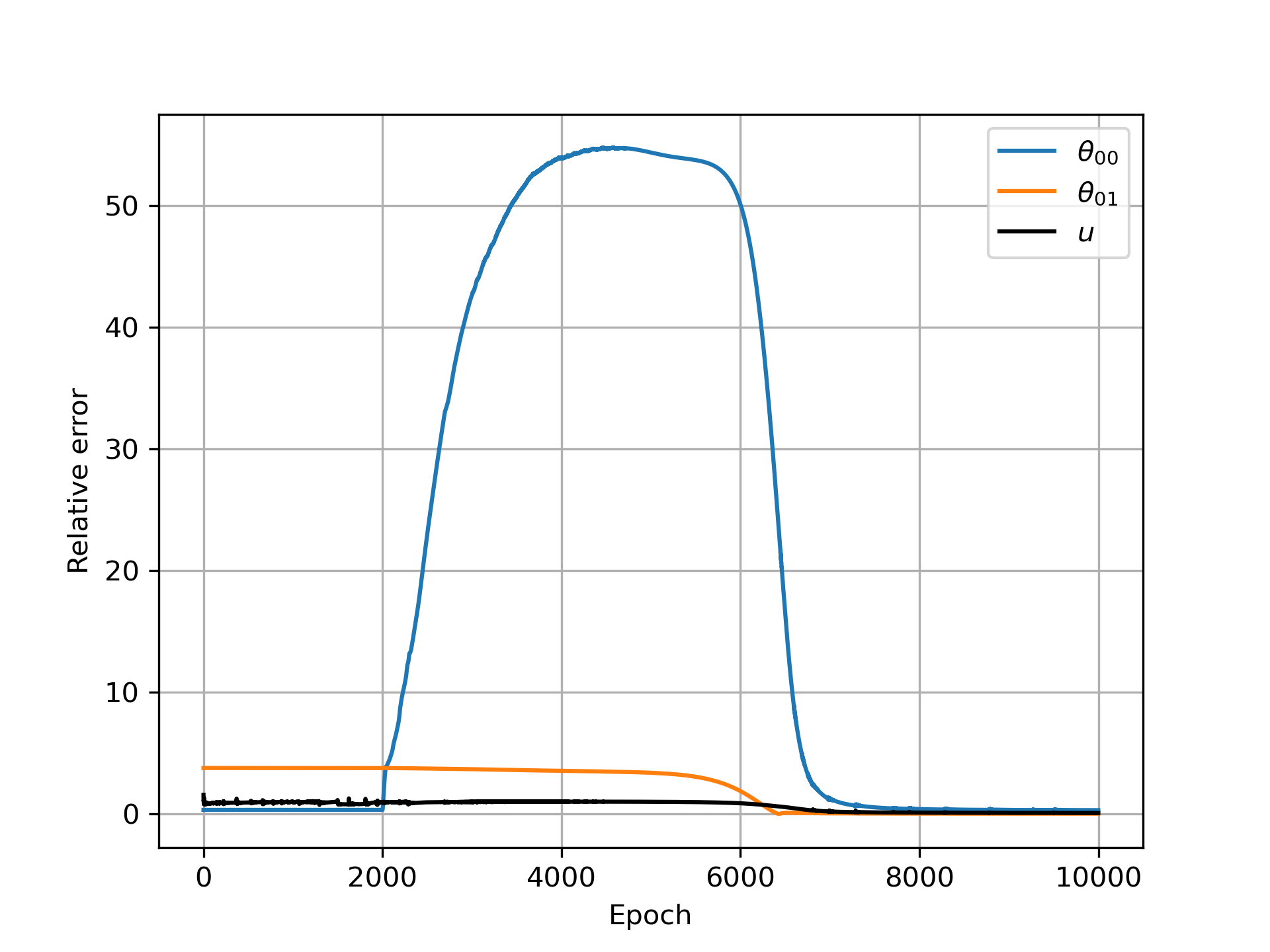}
\includegraphics[width=0.49\textwidth]{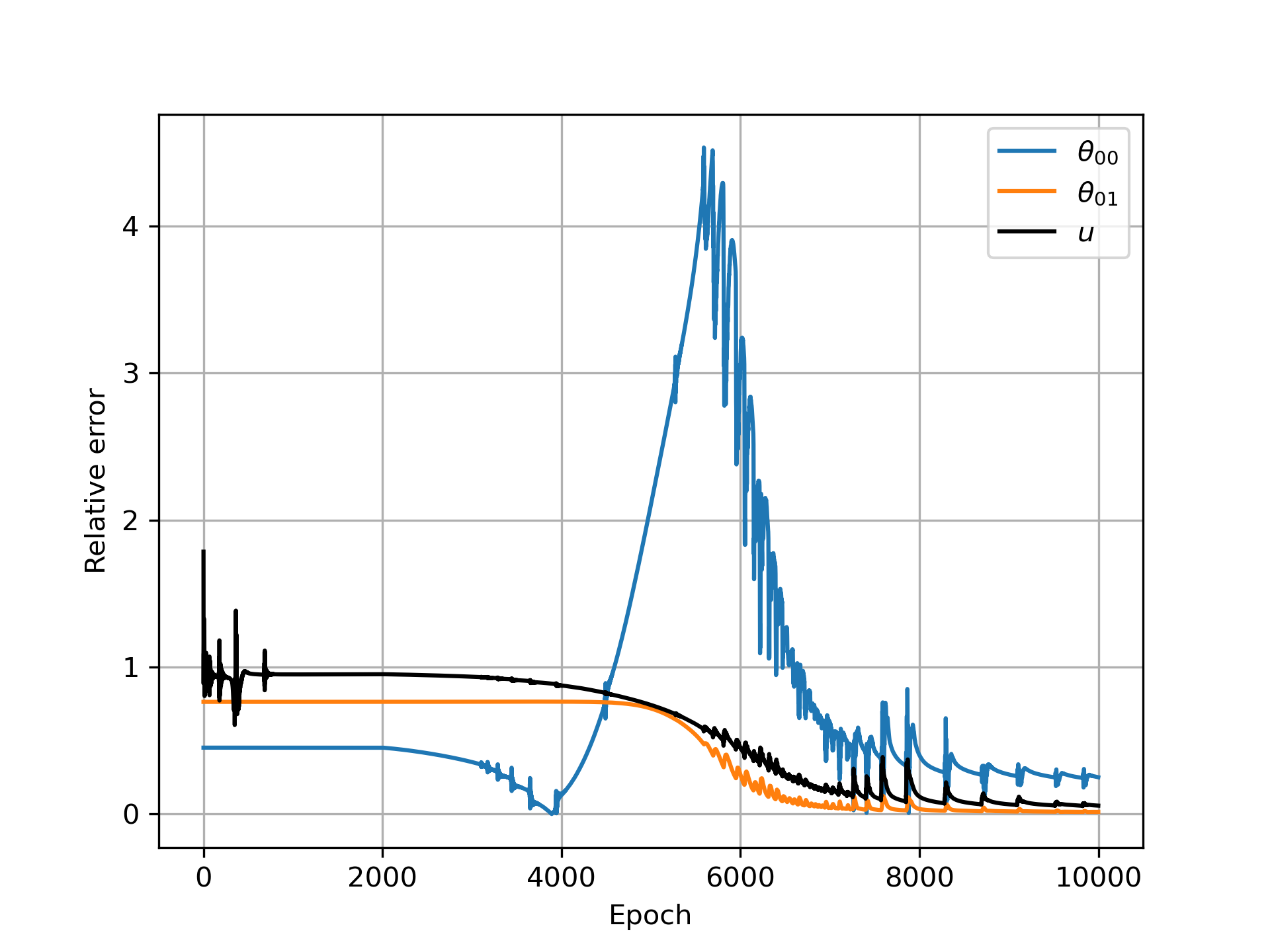}
\caption{Advection-diffusion testcase. Relative error for the advection velocity $\theta_{00}$, dispersion coefficient $\theta_{01}$ and the solution $u$ during the training for two additional random initial value of the parameters.}
\label{fig:appendix_param5}
\end{figure}

\begin{figure}[htbp]
\centering
\includegraphics[width=0.49\textwidth]{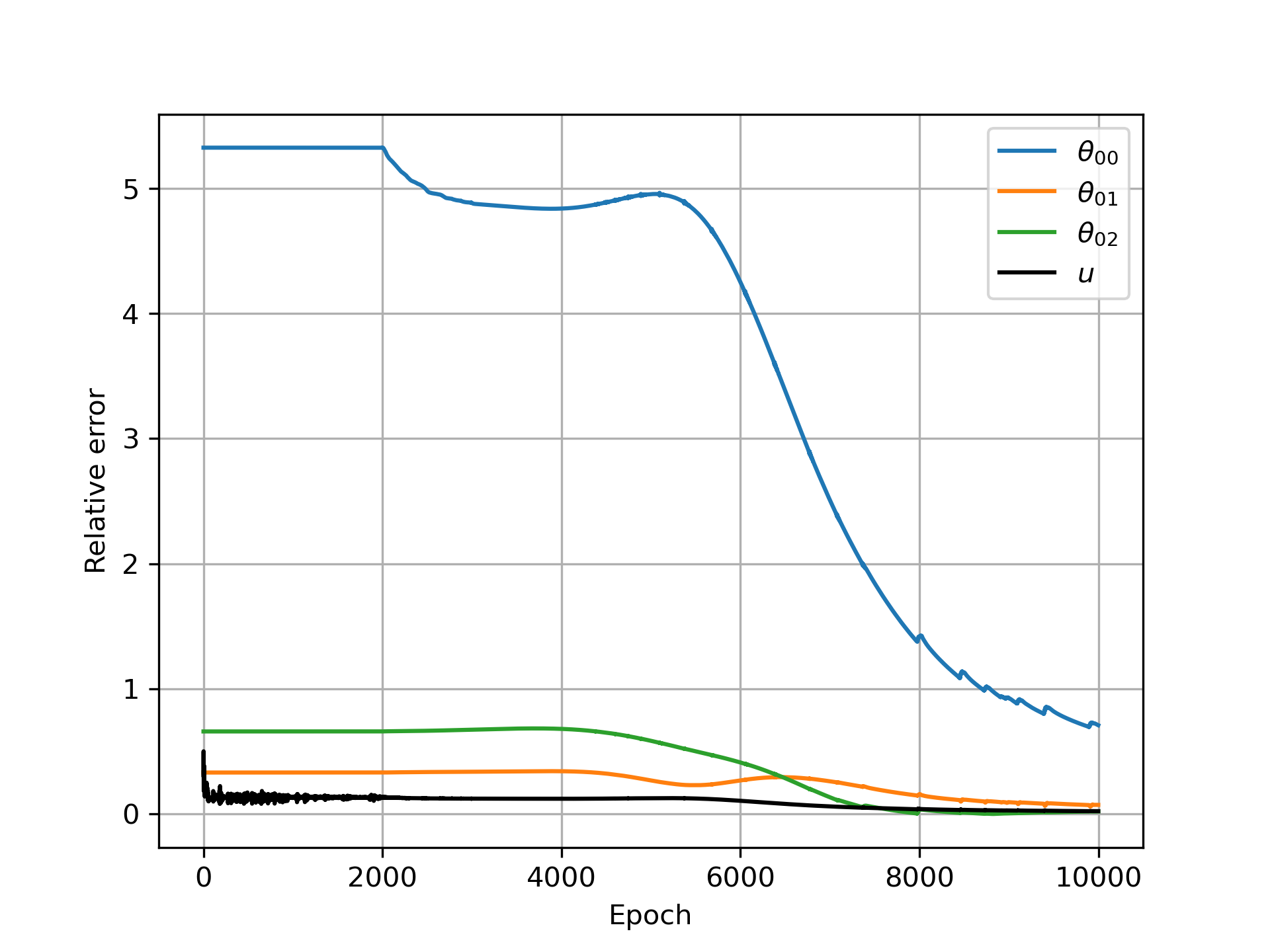}
\includegraphics[width=0.49\textwidth]{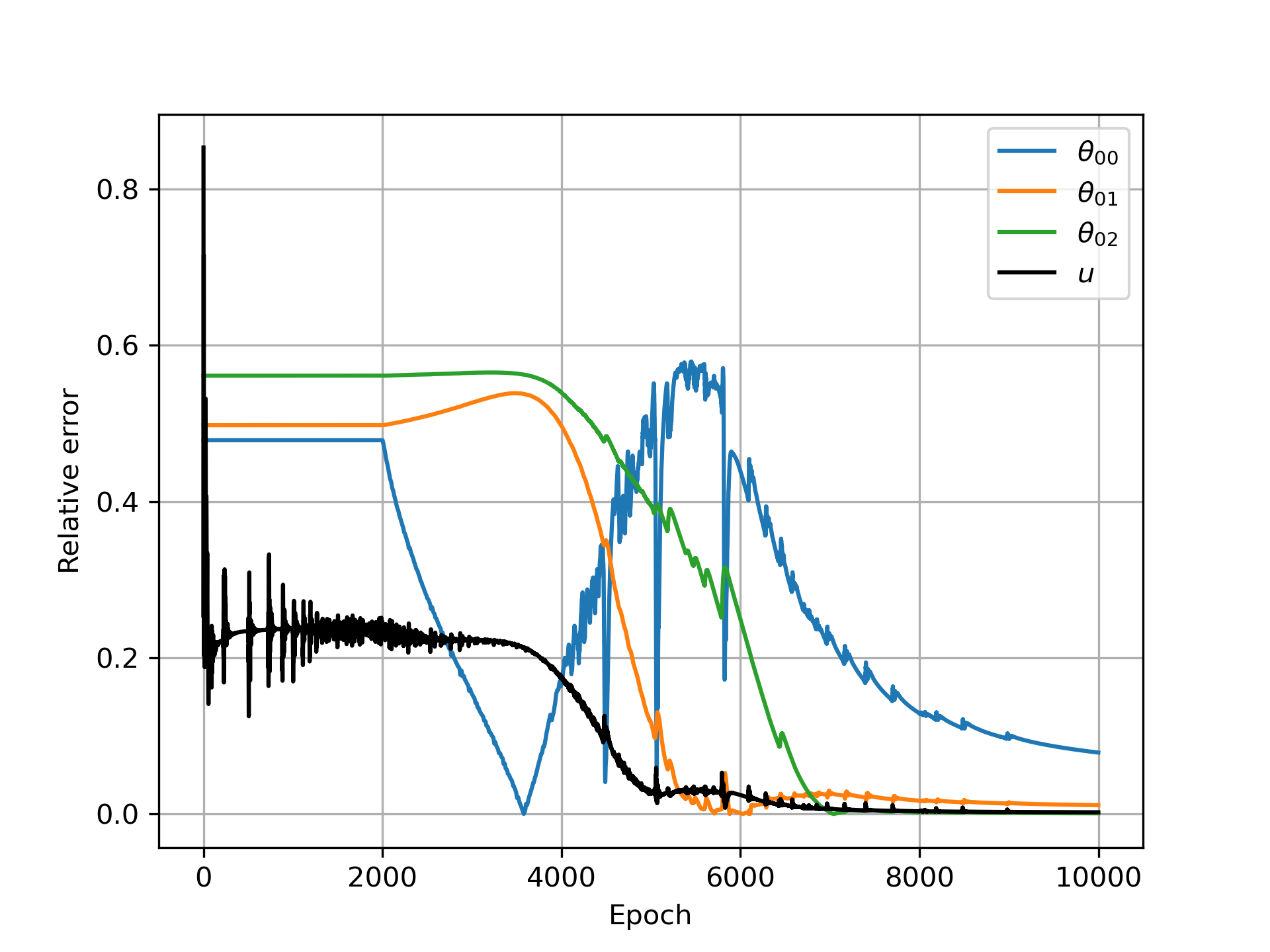}
\caption{Mobile-immobile testcase. Relative error for the dispersion coefficient $\theta_{00}$, advection velocity $\theta_{01}$, and transfer coefficient $\theta_{02}$ and the solution $u$ during the training for two additional random initial value of the parameters.}
\label{fig:appendix_param1}
\end{figure}

Here we present two more random initial conditions for the parameters, for each testcase presented in \cref{sec:results}. The results are shown in \cref{fig:appendix_param0}, \cref{fig:appendix_param5}, \cref{fig:appendix_param1}. Different random seeds are used to initialise the random number generator and parameters are initialised between 0.1 and 10 times the exact values. The PINN is trained with the same hyperparameters as in the main text. The results show that the PINN is able to converge to the correct values of the parameters for different initial conditions, demonstrating the robustness of the proposed architecture.

As it can be seen in \cref{fig:appendix_param0}, \cref{fig:appendix_param5}, \cref{fig:appendix_param1}, the PINN is able to converge to the correct values of the parameters for different initial conditions demonstrating the robustness of the proposed architecture. The interested reader is referred to the code repository for further details on the implementation of the PINN and the training algorithm, and for a more extensive sensitivity analysis with respect to the initial parameters.


\end{document}